%% file: main.tex
\documentclass[accepted]{uai2021} 

\usepackage[british]{babel}

\usepackage{natbib} 
\usepackage{mathtools} 
\usepackage{booktabs} 
\usepackage{tikz} 

\usepackage{amsfonts}
\usepackage{amsmath}
\usetikzlibrary{positioning, arrows}
\usepackage{subcaption}
\usepackage{bbm}
\usepackage{enumitem}
\usepackage[bottom]{footmisc}

\newtheorem{theorem}{Theorem}[section]
\newtheorem{definition}{Definition}[section]

\newtheorem{proposition}{Proposition}[section]

\newcommand{\PP}{\mathbb{P}} 
\newcommand{\EE}{\mathbb{E}} 
\newcommand{\RR}{\mathbb{R}} 
\newcommand{\NN}{\mathbb{N}} 

\newcommand{\indep}{\raisebox{0.05em}{\rotatebox[origin=c]{90}{$\models$}\hspace{0.1em}}}
\newcommand{\dep}{\hspace{0.1em}\raisebox{0.05em}{$\not$}\kern -0.12em \indep}

\newcommand*{\vcenteredhbox}[1]{\begingroup
\setbox0=\hbox{#1}\parbox{\wd0}{\box0}\endgroup}

\title{A Bayesian Nonparametric Conditional Two-sample Test with an Application to Local Causal Discovery}

\author[1]{\href{mailto:philip.boeken@gmail.com?subject=Your UAI 2021 paper}{Philip~A.~Boeken}}
\author[1]{\href{mailto:j.m.mooij@uva.nl}{Joris~M.~Mooij}}

\affil[1]{%
    Korteweg-de Vries Institute for Mathematics\\
    University of Amsterdam
}

\begin{document}
\maketitle

\begin{abstract}
    For a continuous random variable $Z$, testing conditional independence $X\indep Y|Z$ is known to be a particularly hard problem. It constitutes a key ingredient of many constraint-based causal discovery algorithms. These algorithms are often applied to datasets containing binary variables, which indicate the ‘context’ of the observations, e.g.\ a control or treatment group within an experiment. In these settings, conditional independence testing with $X$ or $Y$ binary (and the other continuous) is paramount to the performance of the causal discovery algorithm. To our knowledge no nonparametric ‘mixed’ conditional independence test currently exists, and in practice tests that assume all variables to be continuous are used instead. In this paper we aim to fill this gap, as we combine elements of \citet{holmes2015} and \citet{teymur2020} to propose a novel Bayesian nonparametric conditional two-sample test. Applied to the Local Causal Discovery algorithm, we investigate its performance on both synthetic and real-world data, and compare with state-of-the-art conditional independence tests.
\end{abstract}

\section{Introduction}
\input{input/1_introduction}

\section{Independence testing using P\'olya tree priors}
\input{input/2_polya_trees}

\section{Experiments}
\input{input/3_experiments}

\section{Sensitivity analysis}
\label{section:sensitivity}
\input{input/4_sensitivity_analysis}

\section{Discussion}
\input{input/5_discussion}

\begin{acknowledgements} 
    JMM was supported by the European Research Council (ERC) under the European Union's Horizon 2020 research and
    innovation programme (grant agreement 639466).
\end{acknowledgements}

\bibliography{main}

\newpage
\appendix
\onecolumn
\thispagestyle{empty}
\title{A Bayesian Nonparametric Conditional Two-sample Test with an Application to Local Causal Discovery (Supplementary material)}
\maketitle

\renewcommand{\thesection}{\arabic{section}}
\setcounter{equation}{0}
\setcounter{figure}{0}

\section{Hypothesis testing with P\'olya tree priors}
\input{input/A_polyatrees}

\section{Protein expression data}
\input{input/B_sachs}

\end{document}

%% file: input/1_introduction.tex
Conditional independence testing is a fundamental ingredient of many causal inference algorithms such as the PC algorithm \citep{spirtes1993}, FCI \citep{sprites2013}, and the Local Causal Discovery algorithm \citep{cooper1997}. These algorithms can be proven to be complete, sound, or have other desired properties, but these proofs often invoke the use of an ‘oracle’ for determining conditional independence between variables. In practice, the applicability and performance of the algorithm heavily relies on the reliability of the conditional independence test that is being used. Consequently, incorporating any prior knowledge of the variables involved into the choice of conditional independence test can be desirable.

One way of incorporating prior knowledge is by tailoring the conditional independence tests for $X\indep Y|Z$ on whether the variables involved are discrete or continuous. In the case that the conditioning variable $Z$ is continuous, conditional independence testing is known to be a particularly hard problem \citep{shah2020} and further specifying whether $X$ and $Y$ are continuous or discrete can be beneficial. For the parametric setting multiple ‘mixed’ tests are available \citep{scutari2010, andrews2018, sedgewick2019}. For the nonparametric setting, recent literature proposes multiple tests where $X$ and $Y$ are both assumed to be discrete or both continuous, but to our knowledge no nonparametric test for either $X$ or $Y$ discrete (and the other continuous) currently exists.

Such a ‘mixed’ conditional independence test has a particularly important role in constraint-based causal discovery algorithms that are applied to datasets which are formed by merging datasets from different contexts \citep{mooij2020}. Such a context may for example be whether certain chemicals have been added to a system of proteins (as in Section \ref{section:sachs}), or may be the country of residence of a respondent in an international survey. When certain features of interest (\textit{system variables}) have been measured in different contexts, these measurements can be gathered into a single dataset by adding one or several (often discrete) \textit{context variables} to the dataset, encoding the context that the observation originates from. Merging datasets in this manner may render certain causal relations identifiable, and may improve the reliability of the conditional independence tests due to an increasing sample size \citep{mooij2020}.

Among the continuous conditional independence tests is a recently proposed Bayesian nonparametric test by \citet{teymur2020} which extends a continuous marginal independence test \citep{filippi2017} by utilising \textit{conditional optional P\'olya tree priors} \citep{ma2017}. Although this conditional independence test performs well on data originating from continuous distributions, the prior is misspecified in the case of combinations of discrete and continuous variables. Subsequently, the test has close to zero recall when applied to certain datasets consisting of combinations of discrete and continuous variables.

In this paper we focus on the simplified case of testing $X\indep Y|Z$, where $Z$ and either $X$ or $Y$ is continuous, and the other is binary. We propose a Bayesian nonparametric conditional two-sample test by combining elements of the two-sample test by \citet{holmes2015} and the continuous conditional independence test by \citet{teymur2020}. The two-sample test \citep{holmes2015}, independence test \citep{filippi2017} and our novel conditional two-sample test are empirically compared to both classical and state-of-the-art frequentist (conditional) independence tests when testing for a single (conditional) independence, and when simultaneously testing for multiple (conditional) independences as required by the constraint-based causal discovery algorithm \textit{Local Causal Discovery} (LCD) \citep{cooper1997}.\footnote{Code for the (conditional) independence tests, simulations and results on real world data is publicly available at \url{https://github.com/philipboeken/PTTests}.} Since p-values do not, unlike Bayes factors, reflect any evidence in favour of the null hypothesis, the comparison of Bayesian and frequentist tests in the LCD setting is not straightforward. We propose a measure which allows comparison of the LCD algorithm when using tests from both paradigms, and use it for the comparison of the ensemble of P\'olya tree tests with frequentist tests. We observe that LCD with the ensemble of P\'olya tree tests outperforms other state-of-the-art (conditional) independence tests, while computation time is substantially lower compared to the competing tests.

We apply the LCD algorithm with the P\'olya tree tests to protein expression data from \citet{sachs2005}, and conclude that this implementation provides a result that is more likely to resemble the true model than the output of LCD with the often used partial correlation test.


%% file: input/2_polya_trees.tex
If we let $X:\Omega\rightarrow\mathcal{X}$ be a random variable with distribution $P$ and let $\mathcal{M}$ be the space of all probability distributions on $\mathcal{X}$, then for subsets $\mathcal{M}_0 \subset\mathcal{M}$ and $\mathcal{M}_1\subset\mathcal{M}$ we may test the hypotheses $H_0: P\in\mathcal{M}_0$ and $H_1:P\in\mathcal{M}_1$ by considering \textit{random measures} $\mathcal{P}_0$ and $\mathcal{P}_1$ with distributions $\Pi_0$ and $\Pi_1$ such that $\mathcal{P}_0\in\mathcal{M}_0$ $\Pi_0$-a.s.\ and $\mathcal{P}_1\in\mathcal{M}_1$ $\Pi_1$-a.s. If the posterior distribution of either $\mathcal{P}_0$ or $\mathcal{P}_1$ is consistent (depending on whether $P\in\mathcal{M}_0$ or $P\in\mathcal{M}_1$) and both models $\mathcal{M}_0$ and $\mathcal{M}_1$ are absolutely continuous with respect to some dominating measure, then we may equivalently state the hypotheses as $H_0: X\sim \mathcal{P}_0$ and $H_1:X\sim \mathcal{P}_1$, and test these hypotheses by computing the Bayes factor
\begin{equation}
    \textrm{BF}_{01} = \frac{\PP(H_0)}{\PP(H_1)}\frac{\int_{\mathcal{M}}\prod_{i=1}^n p(X_i)d\Pi_{0}(P)}{\int_{\mathcal{M}}\prod_{i=1}^n p(X_i)d\Pi_{1}(P)},
\end{equation}
where $\PP(H_j)$ is the prior probability of hypothesis $H_j$, $p$ is the Radon-Nikodym derivative of $P$ with respect to the dominating measure, and the integral $\int_{\mathcal{M}}\prod_{i=1}^n p(X_i)d\Pi_{j}(P)$ is the marginal likelihood of the sample $X_1, ..., X_n$ with respect to hypothesis $H_j$. In this work we will use the \textit{P\'olya tree} as a random measure which, under certain assumptions, has a closed form expression for the marginal likelihood of a sample of observations. This is a major benefit compared to e.g.\ the Dirichlet process, as the Dirichlet process often requires costly MCMC sampling to calculate the marginal likelihood.

\begin{figure}
    \centering
    \fbox{\input{graphics/1d_partition}}
    \caption{Construction of a one-dimensional P\'olya tree based on canonical partitions.}
    \label{1d_partition_diagram}
\end{figure}
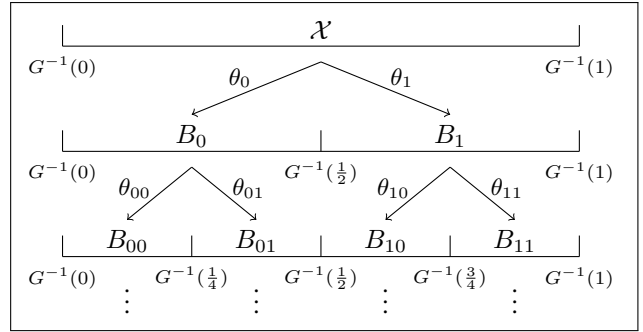

To construct a P\'olya tree on $\mathcal{X} \subseteq \RR$ we consider the set of \textit{canonical partitions of $\mathcal{X}$}, which is defined as the recursive set of partitions
\begin{equation}
    \mathcal{T} = \{\mathcal{X}, \{B_0, B_1\}, \{B_{00}, B_{01}, B_{10}, B_{11}\}, ...\}
\end{equation}
formed by mapping the family of dyadic partitions of $[0,1]$ through the inverse of a cumulative distribution function $G : \mathcal{X}\rightarrow[0,1]$ \citep{ghosal2017}. This results in a family of partitions of $\mathcal{X}$, where for level $j$ we have $\mathcal{X} = \bigcup_{\kappa \in \{0,1\}^j} B_\kappa$, with
\begin{equation}\label{1d_partition}
    B_\kappa := \big[G^{-1}(\tfrac{k-1}{2^j}), G^{-1}(\tfrac{k}{2^j})\big),
\end{equation}
and $k$ denoting the natural number corresponding with the bit string $\kappa \in \{0,1\}^j$. A schematic depiction of this binary tree of partitions is shown in Figure \ref{1d_partition_diagram}. If we define the index set $K := \{\{0,1\}^j: j\in \NN\}$, then the random measure $\mathcal{P}$ is constructed by letting $\mathcal{P}(\mathcal{X}):= 1$ and recursively assigning random probabilities to $B_\kappa\in \mathcal{T}$ by splitting from the mass that is assigned to $B_\kappa$ a fraction $\theta_{\kappa 0}$ to $B_{\kappa 0}$ and a fraction $\theta_{\kappa 1}$ to $B_{\kappa 1}$, where we let $(\theta_{\kappa 0}, \theta_{\kappa 1}) \sim \textrm{Dir}(\alpha_{\kappa 0}, \alpha_{\kappa 1})$. This construction yields a random Borel measure $\mathcal{P}$ on $\mathcal{X}$ \citep{ghosal2017} which adheres to the following definition:

\begin{definition}[\citealp{lavine1992}]
A random probability measure $\mathcal{P}$ on $(\mathcal{X}, \mathcal{B}(\mathcal{X}))$ is said to have a {\normalfont P\'olya tree} distribution with parameter $(\mathcal{T}, \mathcal{A})$, written $\mathcal{P} \sim {\normalfont\textrm{PT}}(\mathcal{T}, \mathcal{A})$, if there exist nonnegative numbers $\mathcal{A} = \{(\alpha_{\kappa 0}, \alpha_{\kappa 1}): \kappa \in K\}$ and random variables $\Theta = \{(\theta_{\kappa 0}, \theta_{\kappa 1}) : \kappa \in K\}$ such that the following hold:
\begin{enumerate}[itemsep=1.8pt, parsep=1.8pt]
    \item all the random variables in $\Theta$ are independent;
    \item for every $\kappa \in K$, we have $(\theta_{\kappa 0}, \theta_{\kappa 1}) \sim {\normalfont\textrm{Dir}}(\alpha_{\kappa 0}, \alpha_{\kappa 1})$;
    \item for every $j \in \NN$ and every $\kappa \in \{0,1\}^j$ we have $\mathcal{P}(B_{\kappa}) = \prod_{i=1}^j \theta_{\kappa_1 ...\kappa_{i}}$.
\end{enumerate}
\end{definition}

Let $X$ be a continuous random variable and consider the P\'olya tree $\mathcal{P} \sim \textrm{PT}(\mathcal{T}, \mathcal{A})$. Drawing a distribution from $\mathcal{P}$ is done by drawing from each of the random variables in $\Theta$. If we let $X_1, ..., X_n$ be a sample from $X$, then the likelihood of that sample with respect to a sampled distribution $\Theta$ from the P\'olya tree $\textrm{PT}(\mathcal{T}, \mathcal{A})$ is
\begin{equation}
    p (X_{1:n}|\Theta, \mathcal{T}, \mathcal{A}) = \prod_{\kappa \in K}\theta_{\kappa}^{n_{\kappa}},
\end{equation}
where $n_\kappa$ denotes the number of observations lying in $B_\kappa$, i.e.~$n_\kappa := |X_{1:n} \cap B_\kappa|$. If we integrate out $\Theta$ we obtain the marginal likelihood
\begin{equation}\label{1d_likelihood}
    p(X_{1:n} | \mathcal{T}, \mathcal{A}) = \prod_{\kappa \in K}\frac{\textrm{B}(\alpha_{\kappa 0} + n_{\kappa 0}, \alpha_{\kappa 1} + n_{\kappa 1})}{\textrm{B}(\alpha_{\kappa 0}, \alpha_{\kappa 1})},
\end{equation}
where $\textrm{B}(\cdot)$ denotes the Beta function.

The choice of $\mathcal{T}$ and $\mathcal{A}$ influences certain characteristics of samples from the P\'olya tree. For example, if we let $\alpha_{\kappa 0} = \alpha_{\kappa 1}$ for all $\kappa \in K$ then the P\'olya tree is centred on the base distribution with cumulative distribution function $G$, i.e.\ $\EE[\mathcal{P}(B_\kappa)] = \int_{B_\kappa} G'(x)dx$. \citet{kraft1964} provides sufficient conditions on $\mathcal{A}$ for the P\'olya tree to be dominated by Lebesgue measure. These conditions are satisfied if for each $\kappa\in \{0,1\}^j$ we take $\alpha_{\kappa } = |\kappa|^2$ with $|\kappa| := j$. The choice of the parameter $\mathcal{A}$ is analysed in Section \ref{section:sensitivity}.

\subsection{A nonparametric conditional two-sample test}
\label{section:cond_two_sample_test}
We now propose a conditional independence test of the type $C\indep X|Z$, where $X$ and $Z$ are continuous one-dimensional random variables and $C$ is a binary random variable. Let $F$ be the conditional distribution of $X|Z$, and let the conditional distributions of $X|\{C = 0\}, Z$ and $X|\{C = 1\}, Z$ be $F^{(0)}$ and $F^{(1)}$ respectively. Then we formulate the conditional independence test between $C$ and $X$ given $Z$ as a two-sample test, i.e.
\begin{align}
    \begin{split}
        & H_0: C \indep X | Z \iff F^{(0)} = F^{(1)} = F \\
        & H_1: C \dep X | Z \iff F^{(0)} \neq F^{(1)}.
    \end{split}
\end{align}

Following \citet{teymur2020} we will utilise the \textit{conditional optional P\'olya tree} (cond-OPT) prior \citep{ma2017} for modelling the conditional distributions $F$, $F^{(0)}$ and $F^{(1)}$. The cond-OPT is a random conditional probability measure on e.g.\ $\mathcal{X}\times \mathcal{Z}$, where $X$ is the response variable and $Z$ is the predictor. In order to construct the cond-OPT, we first construct a family of partitions $\mathcal{T}_Z$ of $\mathcal{Z}$ according to the partitioning scheme of the optional P\'olya tree (OPT) \citep{wong2010}, which results in a random subset of the canonical partitions $\mathcal{T}$ as constructed by equation (\ref{1d_partition}). This random subset of $\mathcal{T}$ is obtained by first adding $B_\emptyset := \mathcal{Z}$ to $\mathcal{T}_Z$. Then we sample from the random variable $S \sim \textrm{Bernoulli}(\rho)$; if $S=1$ we stop the partitioning procedure, and if $S=0$ we add $B_0$ and $B_1$ to $\mathcal{T}_Z$. Then, for both $B_0$ and $B_1$ we repeat this procedure; we first draw $S$ from $\textrm{Bernoulli}(\rho)$ and depending on the outcome we add the children of $B_{0}$, then we repeat this to possibly add the children of $B_{1}$. This process is iterated, and terminates a.s. when $\rho > 0$.

Having obtained the family $\mathcal{T}_Z$, we construct a ‘local’ random measure $\mathcal{P}(\cdot |B_\kappa)$ on $\mathcal{X}$ for each $B_\kappa\in\mathcal{T}_Z$ by letting $\mathcal{P}(\cdot|B_\kappa)\sim\textrm{PT}(\mathcal{T}, \mathcal{A})$, and we define the conditional probability $\mathcal{P}(\cdot | Z=z)$ to be constant and equal to the local P\'olya tree $\mathcal{P}(\cdot|B_\kappa)$ on the stopped set $B_\kappa \ni z$. The resulting family of random measures on $\mathcal{X}$ is the conditional optional P\'olya tree (cond-OPT) \citep{ma2017}. When using the canonical partitions for both $\mathcal{X}$ and $\mathcal{Z}$ and assuming that all the local P\'olya trees are a.s.\ dominated by Lebesgue measure, \citet{ma2017} shows that the cond-OPT places positive probability on all $L_1$ neighbourhoods of any conditional density $f(\cdot|\cdot)$ on $\mathcal{X} \times \mathcal{Z}$.

When we are given $n$ i.i.d.\ observations $(C_1, X_1, Z_1), ..., (C_n, X_n, Z_n)$, then under the null hypothesis we are interested in the marginal likelihood of a sample $(X_1, Z_1), ..., (X_n, Z_n)$ with respect to the cond-OPT prior. This is obtained by for every $B_\kappa\in \mathcal{T}_Z$ considering the subsample $X(B_\kappa) := \{X_j : Z_j \in B_\kappa\}$. As the cond-OPT prior considers a general P\'olya tree prior for this subsample, we simply compute the marginal likelihood 
\begin{equation}
    p_X(B_\kappa) := p(X(B_\kappa) | \mathcal{T}, \mathcal{A})
\end{equation}
using equation (\ref{1d_likelihood}). If $B_\kappa$ is a so called leaf-set, i.e.\ the set contains at most one observation or it has no children in the family of partitions $\mathcal{T}_Z$, then we simply return this marginal likelihood. If $B_\kappa$ is not a leaf-set, we continue along the children $B_{\kappa 0}$ and $B_{\kappa 1}$. We integrate out the randomness of the random family of partitions by considering the entire family of canonical partitions $\mathcal{T}$ of $\mathcal{Z}$, and incorporating the stopping probabilities $S$ by weighing the elements $B_\kappa \in\mathcal{T}$ of level $|\kappa|$ with $\EE(1-S)^{|\kappa|} = (1-\rho)^{|\kappa|}$. The recursive mixing formula is given by
\begin{equation*}\begin{split}
  \Phi_X(B_\kappa) = \begin{cases}
      p_X(B_\kappa) \ \qquad\qquad\qquad\text{if $B_\kappa$ is a leaf-set} \\
    \begin{split}
      &\rho \cdot p_X(B_\kappa) + {} \\
      &(1-\rho) \cdot \Phi_X(B_{\kappa 0})\Phi_X(B_{\kappa 1})
    \end{split} \quad \text{otherwise,}
    \end{cases}
\end{split}\end{equation*}
and the resulting quantity $\Phi_X(B_\kappa)$ is the marginal likelihood of $\{(X_1, Z_1), ..., (X_n, Z_n)\}\cap \mathcal{X} \times  B_\kappa$, with respect to the cond-OPT. 

Under the alternative hypothesis we split the sample into sets $X^{(0)}:= \{(X_j, Z_j) : C_j = 0\}$ and $X^{(1)}:= \{(X_j, Z_j) : C_j = 1\}$, and compute the marginal likelihoods $\Phi_{X^{(0)}}(\mathcal{Z})$ and $\Phi_{X^{(1)}}(\mathcal{Z})$ of these sets with respect to (independent) cond-OPT priors. 
We finally test the hypothesis by computing the Bayes factor
\begin{equation}
    \textrm{BF}_{01} = \frac{\Phi_X(\mathcal{Z})}{\Phi_{X^{(0)}}(\mathcal{Z})\Phi_{X^{(1)}}(\mathcal{Z})},
\end{equation}
where we have set the prior odds to 1.
 
We note that when no data is provided for $Z$ and thus $\mathcal{Z}$ constitutes a leaf-set, this test defaults to the two-sample test from \citet{holmes2015}. An overview of this two-sample test and the continuous independence test by \citet{filippi2017} is provided in the supplement.

%% file: graphics/1d_partition.tex
\begin{tikzpicture}[xscale=1.7, yscale=0.7]

\node () at (0, 2.5) {};

\node (omega) at (0, 2.3) {$\mathcal{X}$};
\draw (-2, 2) -- (2, 2) ;
\draw (-2, 2) -- (-2, 2.4) ;
\draw (2, 2) -- (2, 2.4) ;
\node () at (-2, 1.6) {\scriptsize $G^{-1}(0)$};
\node () at (2, 1.6) {\scriptsize $G^{-1}(1)$};

\draw[->] (0, 1.7) -- (-1, 0.7) node [pos=0.3, left=0.8em] {\footnotesize$\theta_0$};
\draw[->] (0, 1.7) -- (1, 0.7) node [pos=0.3, right=0.8em] {\footnotesize$\theta_1$};

\node (B_0) at (-1, 0.3) {$B_0$};
\node (B_1) at (1, 0.3) {$B_1$};
\draw (-2, 0) -- (2, 0) ;
\draw (-2, 0) -- (-2, 0.4) ;
\draw (0, 0) -- (0, 0.4) ;
\draw (2, 0) -- (2, 0.4) ;
\node () at (-2, -0.4) {\scriptsize $G^{-1}(0)$};
\node () at (0, -0.4) {\scriptsize $G^{-1}(\tfrac{1}{2})$};
\node () at (2, -0.4) {\scriptsize $G^{-1}(1)$};

\draw[->] (-1, -0.3) -- (-1.5, -1.3) node [pos=0.4, left=0.2em] {\footnotesize$\theta_{00}$};
\draw[->] (-1, -0.3) -- (-0.5, -1.3) node [pos=0.4, right=0.2em] {\footnotesize$\theta_{01}$};
\draw[->] (1, -0.3) -- (0.5, -1.3) node [pos=0.4, left=0.2em] {\footnotesize$\theta_{10}$};
\draw[->] (1, -0.3) -- (1.5, -1.3) node [pos=0.4, right=0.2em] {\footnotesize$\theta_{11}$};

\node (B_00) at (-1.5, -1.7) {$B_{00}$};
\node (B_01) at (-0.5, -1.7) {$B_{01}$};
\node (B_10) at (0.5, -1.7) {$B_{10}$};
\node (B_11) at (1.5, -1.7) {$B_{11}$};
\draw (-2, -2) -- (2, -2) ;
\draw (-2, -2) -- (-2, -1.6) ;
\draw (-1, -2) -- (-1, -1.6) ;
\draw (0, -2) -- (0, -1.6) ;
\draw (1, -2) -- (1, -1.6) ;
\draw (2, -2) -- (2, -1.6) ;
\node () at (-2, -2.4) {\scriptsize $G^{-1}(0)$};
\node () at (-1, -2.4) {\scriptsize $G^{-1}(\tfrac{1}{4})$};
\node () at (0, -2.4) {\scriptsize $G^{-1}(\tfrac{1}{2})$};
\node () at (1, -2.4) {\scriptsize $G^{-1}(\tfrac{3}{4})$};
\node () at (2, -2.4) {\scriptsize $G^{-1}(1)$};


\node () at (-1.5, -2.7) {$\vdots$};
\node () at (-0.5, -2.7) {$\vdots$};
\node () at (0.5, -2.7) {$\vdots$};
\node () at (1.5, -2.7) {$\vdots$};

\node () at (0, -3) {};

\end{tikzpicture}

%% file: input/3_experiments.tex
Implementing the conditional independence test requires choosing certain hyperparameters. As mentioned earlier, we set $\alpha_{\kappa} = |\kappa|^2$. As argued by \citet{lavine1994} we will only consider partitions up to a pre-determined level $J$, making $\mathcal{P}$ into a \textit{truncated P\'olya Tree}. \citet{hanson2002} provide the rule of thumb $J = \lfloor\log_2(n)\rfloor$, which corresponds to on average finding one observation in each element of the partition. We find however that $J = \lfloor\log_4(n)\rfloor$, which corresponds to finding approximately $\sqrt{n}$ observations in each element of the partition, provides similar results and considerably reduces computation time, so we use this maximum depth. Throughout this work we will use the standard Gaussian cdf $G$ to form the canonical partitions. In conjunction with this mean measure, we standardise the data before computing the marginal likelihoods. For computing the marginal likelihood of the cond-OPT we use $\rho = 1/2$ \citep{ma2017}. Similar to the computation of marginal likelihoods of regular P\'olya trees, we use a maximum partitioning depth of $\lfloor\log_4(n)\rfloor$, so we consider $B_\kappa \in \mathcal{T}_Z$ to be a leaf-set if it contains at most one value, or if $|\kappa| = \lfloor\log_4(n)\rfloor$.

All experiments are run on a MacBook Pro with a 3.1 GHz CPU and 16GB of RAM, with a parallelised R implementation of the LCD algorithm.
Code for the (conditional) independence tests, simulations and results on real world data is publicly available at \url{https://github.com/philipboeken/PTTests}.

\subsection{Local Causal Discovery}
As mentioned earlier, a ‘mixed’ conditional independence test as proposed in Section \ref{section:cond_two_sample_test} is specifically needed when applying causal discovery algorithms to datasets containing binary (or discrete) \textit{context variables}, which encode the context that observations of the \textit{system variables} (the variables of interest) originate from. In accordance with \citet{mooij2020}, we regard both the context variables and the system variables as distributed according to the solution of a \textit{Structural Causal Model} (SCM) \citep{pearl2009}. A relatively insightful causal discovery algorithm is the \textit{Local Causal Discovery} (LCD) algorithm \citep{cooper1997}. Although often referred to as an algorithm, it essentially consists of the following proposition:
\begin{proposition}[LCD, \cite{mooij2020}]\label{LCD}
    If the data generating process of the triple of random variables $(C, X, Y)$ has no selection bias, can be modelled by a faithful simple SCM, and $X$ is not a cause of $C$, then the presence of (in)dependences
    \begin{equation}
        C\dep X, \quad X \dep Y, \quad C\indep Y|X
    \end{equation}
    implies that $X$ is a (possibly indirect) cause of $Y$. If this is the case, we speak of the `LCD triple' $(C, X, Y)$.
\end{proposition}
By repeatedly applying this proposition to different triples of random variables one can (partially) reconstruct the underlying causal graph of the dataset at hand. If we are provided with a dataset consisting of observations of context variables $(C_k)_{k\in\mathcal{K}}$ for some index set $\mathcal{K}$ and system variables $(X_i)_{i\in\mathcal{I}}$ for some index set $\mathcal{I}$ for which we assume that the system variables do not cause the context variables, then we may iteratively apply Proposition \ref{LCD} to all triples $(C_k, X_i, X_{i'})$ where $k\in\mathcal{K}$ and $i\neq i' \in \mathcal{I}$, and provide a directed graph as output where the edges can be interpreted as representing indirect causal effects.

\begin{figure*}[!htb]
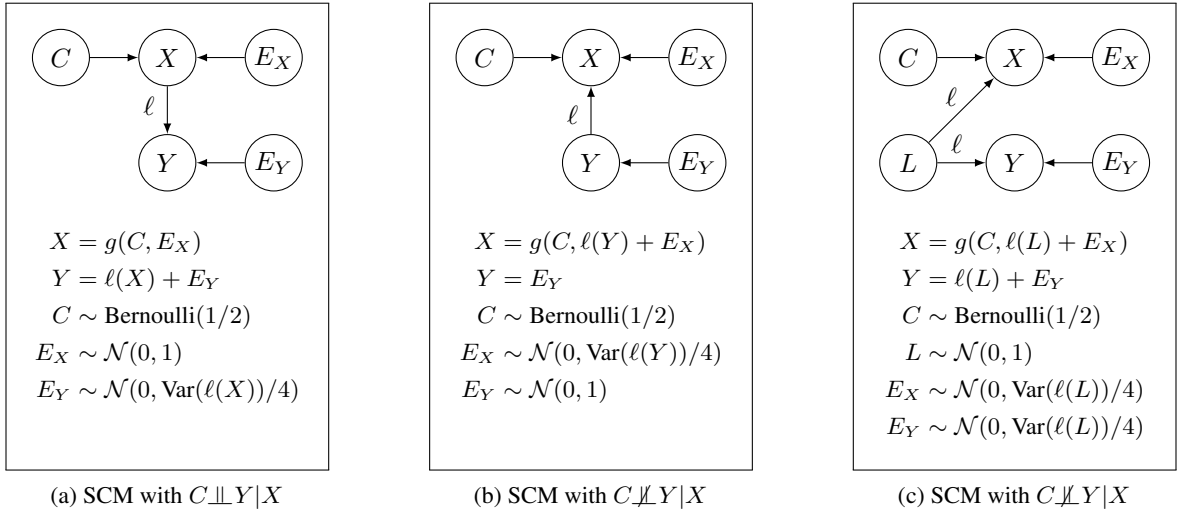

    \centering
    \begin{subfigure}[b]{0.32\textwidth}
        \centering
        \include{graphics/graph1}
        \caption{SCM with $C\indep Y | X$}
        \label{graphs:1}
    \end{subfigure}
    \begin{subfigure}[b]{0.32\textwidth}
        \centering
        \include{graphics/graph2}
        \caption{SCM with $C\dep Y | X$}
        \label{graphs:2}
    \end{subfigure}
    \begin{subfigure}[b]{0.32\textwidth}
        \centering
        \include{graphics/graph3}
        \caption{SCM with $C\dep Y | X$}
        \label{graphs:3}
    \end{subfigure}
    \caption{Three SCMs used for the simulations.}
    \label{fig:graphs}
\end{figure*}

\subsection{Simulations}
\label{simulations}
In our simulations we repeatedly simulate a triple of random variables $(C, X, Y)$. Each time we simulate a set of observations, we test for $C\indep X, X\indep Y$ and $C\indep Y|X$ individually, and by combining the output of these three tests we formulate the output of the LCD algorithm. Upon repeating this scheme a number of times we are able to display ROC curves for each of the three test cases, and for the LCD algorithm. To widen the scope of this setup, in each round of simulations we randomly choose one of the graphs of Figure \ref{fig:graphs}, and we randomly pick the relations between $C$ and $X$ and between $X$ and $Y$ from predefined, varying possibilities. More specifically, if we let $E$ be an external factor (possibly depending on $Y$) we set $X$ equal to $g(C, E)$, which is randomly chosen from
\begin{align}
    g(c, e) =
    \begin{cases}
        e                      & \text{no intervention}      \\
        (1-c)e + c(e + \theta) & \text{mean shift}           \\
        (1-c)e + c(1+\theta)e  & \text{variance shift}       \\
        (1-c)e + c\theta       & \text{perfect intervention} \\
        (1-c)e + c(e + B)      & \text{mean shift mixture},
    \end{cases}
\end{align}
with $\theta \sim \mathcal{U}(\{2, 3, 4, 5, 6\})$ independently drawn per round of simulations and $B \sim \mathcal{U}(\{-1, \theta\})$ independently drawn for every observation. These mappings between $C$ and $X$ can be interpreted as setting $X$ equal to the value $E$ in context $C=0$, and intervening on $X$ in context $C=1$. If we for example inspect the ‘mean shift’, then if $C=1$ we intervene on the distribution of $X$ by shifting the mean of $X$ with the amount $\theta$. When simulating multiple observations, this intervention on $X$ is performed on approximately half of these observations, due to $C$ having a Bernoulli(1/2) distribution. The relation $\ell$ between $X$ and $Y$ is randomly picked from
\begin{align}
    \ell(x) =
    \begin{cases}
        0 \quad                     & \text{no link}    \\
        x \quad                     & \text{linear}     \\
        x^2 \quad                   & \text{parabolic}  \\
        \sin(12\pi \tilde{x}) \quad & \text{sinusoidal}
    \end{cases}
\end{align}
where $\tilde{x} = x / (\max(x_1,..., x_n) - \min(x_1, ... ,x_n))$. It depends on which graph from Figure \ref{fig:graphs} is chosen whether we have $X\overset{\ell}{\rightarrow} Y$, $X\overset{\ell}{\leftarrow} Y$ or $X\overset{\ell}{\leftarrow} L \overset{\ell}{\rightarrow} Y$, where in the last case the two $\ell$'s are drawn independently. The possibility of picking $g(c, e)=e$ and $\ell(x) = 0$ ensures the occurrence of $C\indep X$ and $X\indep Y$ respectively, which in turn enables plotting ROC curves of these test cases.

\begin{figure*}[ht]
    \centering
    \begin{subfigure}[t]{0.245\textwidth}
        \centering
        \includegraphics[width=1\textwidth]{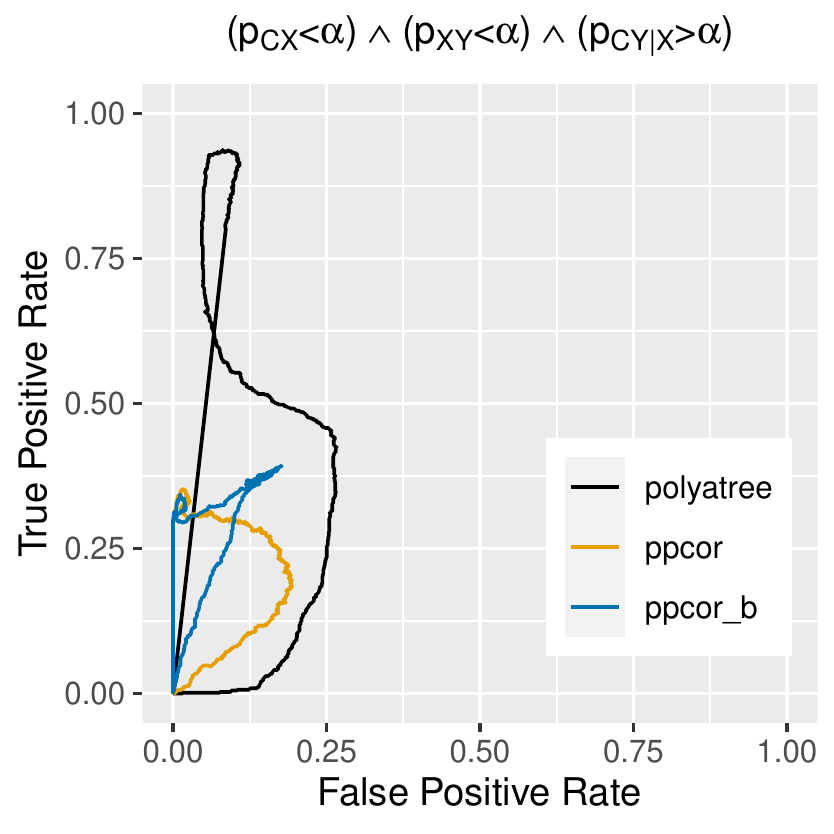}
        \caption{}
        \label{fig:plot1}
    \end{subfigure}
    \begin{subfigure}[t]{0.245\textwidth}
        \centering
        \includegraphics[width=1\textwidth]{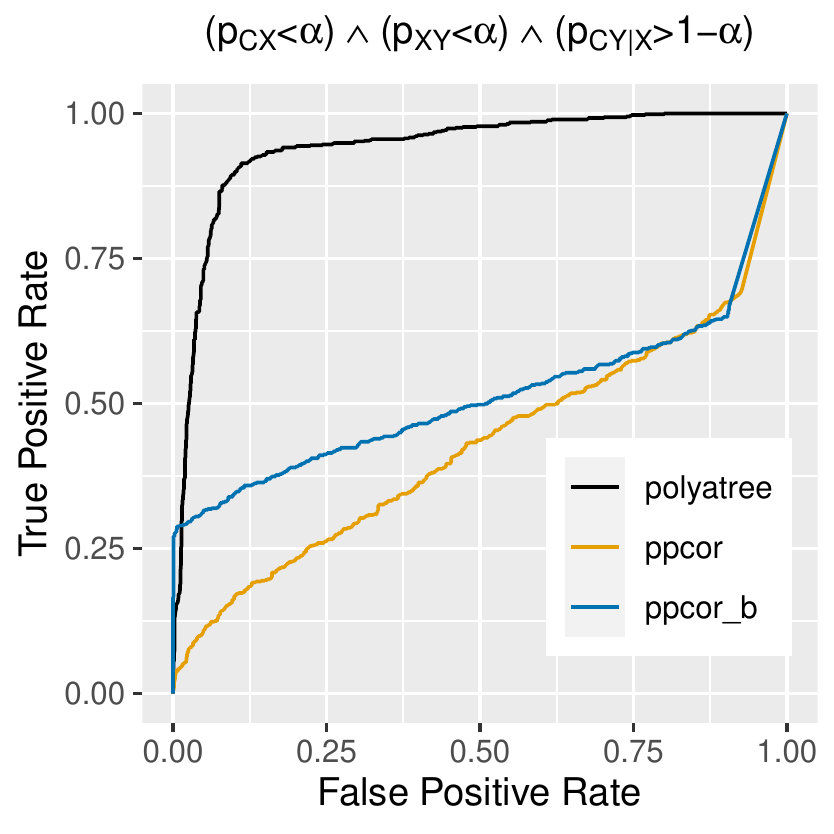}
        \caption{}
        \label{fig:plot2}
    \end{subfigure}
    \begin{subfigure}[t]{0.245\textwidth}
        \centering
        \includegraphics[width=1\textwidth]{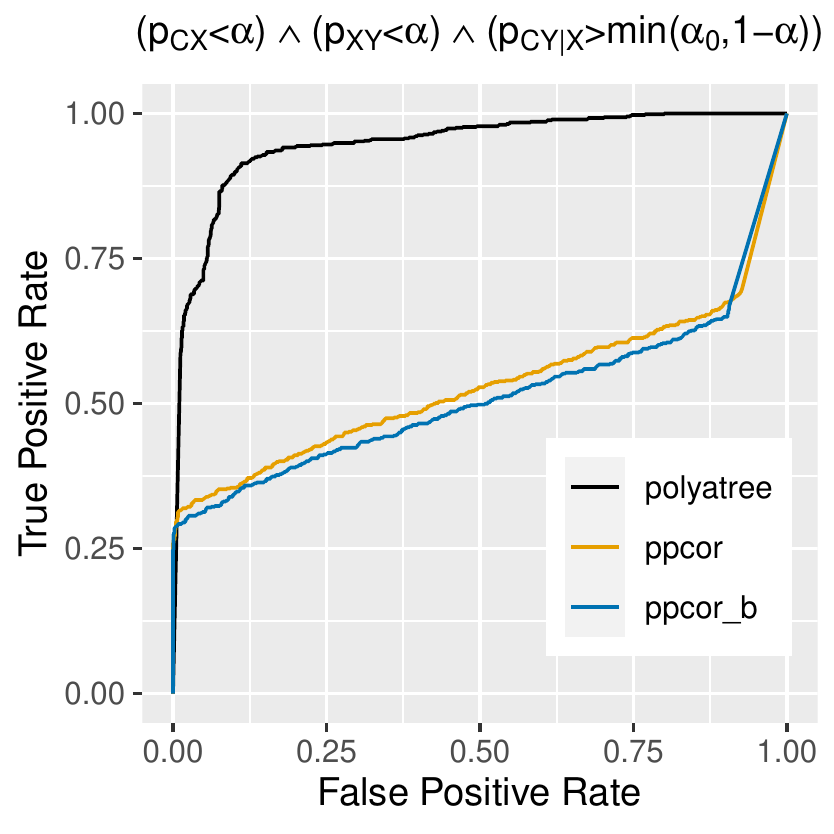}
        \caption{}
        \label{fig:plot0}
    \end{subfigure}
    \caption{ROC curves of different ways of scoring an LCD triple $(C, X, Y)$. See main text for details.}
    \label{lcd_measures}
\end{figure*}

We compare the P\'olya tree based ensemble of the two-sample test \citep{holmes2015}, independence test \citep{filippi2017} and conditional two-sample test (Section \ref{section:cond_two_sample_test}), denoted by \verb+polyatree+, with both classical and recently proposed (conditional) independence tests. The tests that are suitable for mixed testing are \verb+mi_mixed+ and \verb+lr_mixed+, where the former is based on mutual information and uses the implementation of the \verb+bnlearn+ package \cite{scutari2010}, and where the latter is a likelihood ratio test of linear and logistic regressions \citep{sedgewick2019}. Among the more classical continuous tests is the Pearson correlation- and partial correlation test, denoted by \verb+ppcor+, implemented using the synonymous R-package \citep{kim2015}. \citet{harris2013} promote the use of Spearman's (partial) rank correlation test in the context of nonparanormal models, which we denote by \verb+spcor+. Among the more state-of-the-art continuous tests is the \textit{Generalised Covariance Measure} (GCM) \citep{shah2020}, which can be loosely interpreted as a nonlinear extension of the partial correlation test. The GCM is implemented with penalised regression splines as provided by the R-package \verb+GeneralisedCovarianceMeasure+, and is denoted by \verb+gcm+. Departing from the regression-type independence tests, we also consider the \textit{Randomised Conditional Correlation Test} (RCoT) as proposed by \citet{strobl2019}, which closely approximates the Kernel Conditional Independence test by \citet{zhang2011} at the benefit of significantly lower computation time. For marginal independence testing the RCoT defaults to an approximate version of the Hilbert-Schmidt Independence Criterion \citep{gretton2008}. This ensemble is denoted by \verb+rcot+. Lastly we compare to the \textit{Classifier Conditional Independence Test} (CCIT) \citep{sen2017}, denoted by \verb+ccit+, which uses the XGBoost binary classifier to assess presence of conditional independence.

Comparing Bayesian and frequentist tests based on their performance in the LCD algorithm is not straightforward, since the triple of tests for $C\dep X, X\dep Y$ and $C\indep Y|X$ does not by default output a single confidence score. For each test we output the p-value, or in case of the Bayesian tests the $H_0$ model evidence $\PP(H_0|\textrm{data})$.\footnote{Recall that $\PP(H_0|\textrm{data}) = 1 - (1 + \textrm{BF}_{01})^{-1}$.} We construct ROC curves for testing `positive' outcomes $C\dep X$, $X\dep Y$ and $C\dep Y|X$ by varying the threshold $\alpha$, representing the upper bound on the p-value/model evidence for drawing a positive conclusion. The triple $(C, X, Y)$ is given a ‘positive' label if the data is generated according to the relation $C\rightarrow X \rightarrow Y$. Typically, varying the threshold $\alpha$ from 0 to 1 produces an ROC curve between the points $(0,0)$ and $(1,1)$. If we denote the frequentist p-values or Bayesian $H_0$ model evidence for the tests $C\indep X$, $X\indep Y$ and $C\indep Y|X$ with $p_{CX}$, $p_{XY}$ and $p_{CY|X}$ respectively (with independence under the null hypothesis), and if we were to use the same $\alpha$ as threshold for testing whether $p_{CX} < \alpha$, $p_{XY} < \alpha$ and $p_{CY|X} > \alpha$, then varying $\alpha$ between 0 and 1 does not result in a curve between $(0,0)$ and $(1,1)$, as shown in Figure \ref{fig:plot1}.  To assess whether we provide a fair comparison between Bayesian and frequentist tests, we include a Bayesian version of the Pearson (partial) correlation test \citep{wetzels2012}, denoted by \verb+ppcor_b+. Alternatively we could use $\alpha$ for testing $p_{CX} < \alpha$, $p_{XY} < \alpha$ and $p_{CY|X} > 1-\alpha$, as shown in Figure \ref{fig:plot2}. In this case the level $\alpha$ reflects the amount of evidence for the desired conclusions $C\dep X$, $X\dep Y$ and $C\indep Y|X$. For frequentist tests this would not make sense, as for decreasing $\alpha$ we require more evidence for $H_0:C\indep Y|X$, and the p-value has a uniform distribution under $H_0$. This is remedied by, when testing for independence $C\indep Y|X$, only varying $\alpha$ between $0$ and a fixed $\alpha_0$ (Figure \ref{fig:plot0}). More specifically, for level $\alpha$ the LCD algorithm outputs the score
\begin{multline}
    s_{\textrm{LCD}} = \mathbbm{1}_{[0,\alpha]}(p_{CX}) \cdot \mathbbm{1}_{[0,\alpha]}(p_{XY}) \\
    \cdot \mathbbm{1}_{(\alpha_0, 1] \cup (1-\alpha, 1]}(p_{CY|X}),
\end{multline}
where we let $\alpha_0 = 0.05$ for frequentist tests and $\alpha_0=0.5$ for Bayesian tests. Although this $\alpha_0$ is quite arbitrarily chosen, the use of this performance measure is corroborated by the observation that in Figure \ref{fig:plot0} the frequentist partial correlation and Bayesian partial correlation tests have similar performance.

\begin{figure*}[ht]
    \centering
    \begin{subfigure}[t]{0.245\textwidth}
        \centering
        \includegraphics[width=1\textwidth]{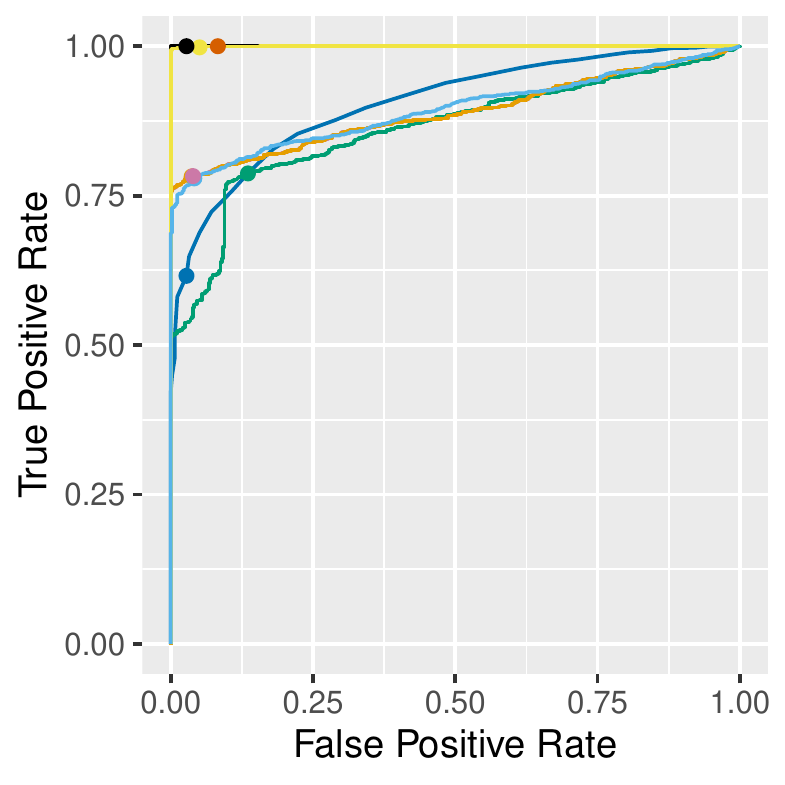}
        \vspace{-5mm}
        \caption{$C \dep X$}
        \label{fig:roc_ts}
    \end{subfigure}
    \begin{subfigure}[t]{0.245\textwidth}
        \centering
        \includegraphics[width=1\textwidth]{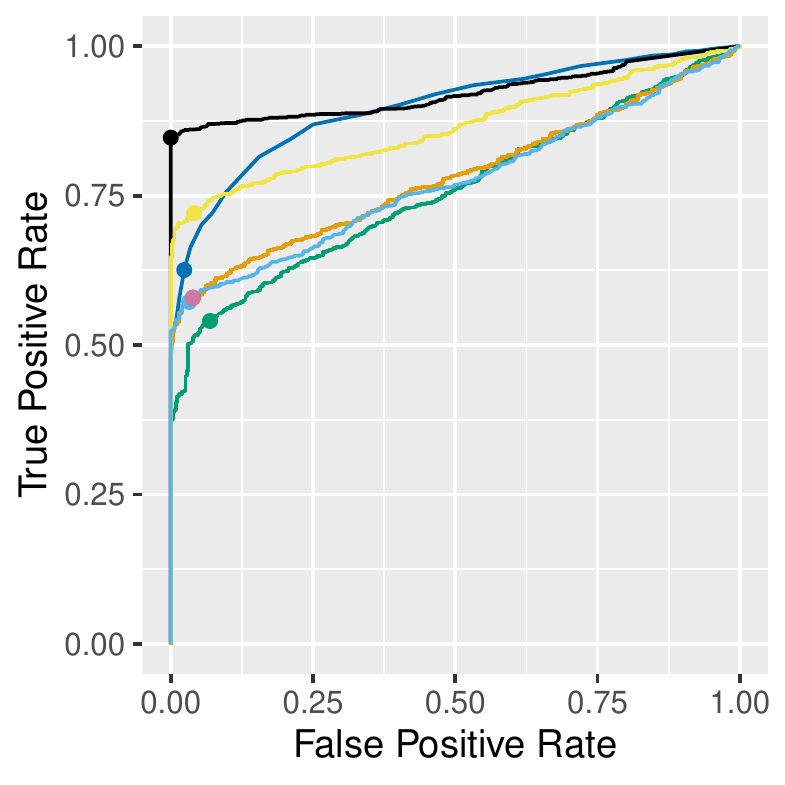}
        \vspace{-5mm}
        \caption{$X \dep Y$}
        \label{fig:roc_uci}
    \end{subfigure}
    \begin{subfigure}[t]{0.245\textwidth}
        \centering
        \includegraphics[width=1\textwidth]{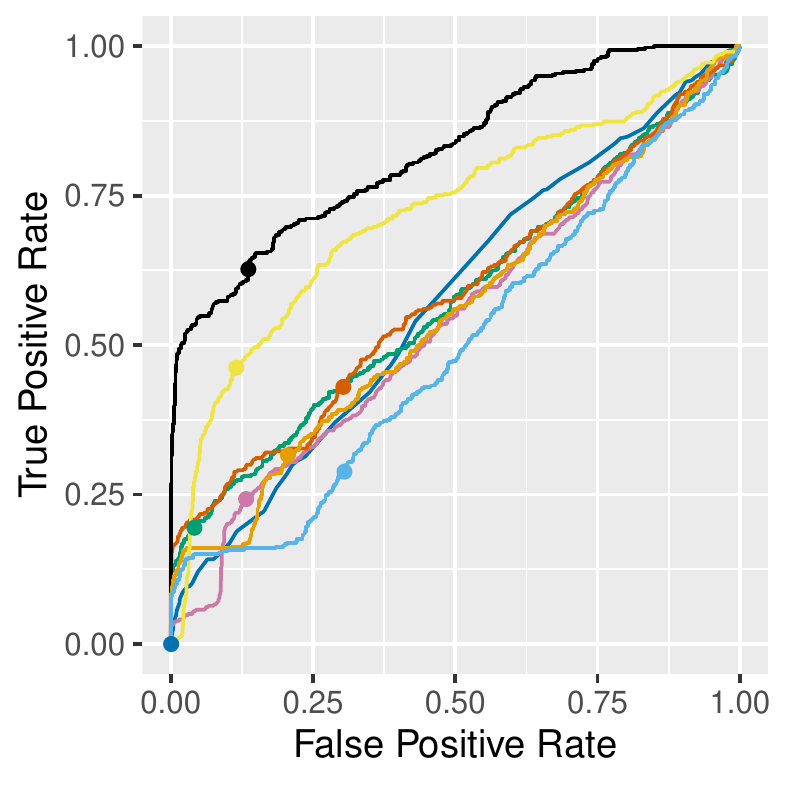}
        \vspace{-5mm}
        \caption{$C \dep Y|X$}
        \label{fig:roc_ci}
    \end{subfigure}
    \begin{subfigure}[t]{0.245\textwidth}
        \centering
        \includegraphics[width=1\textwidth]{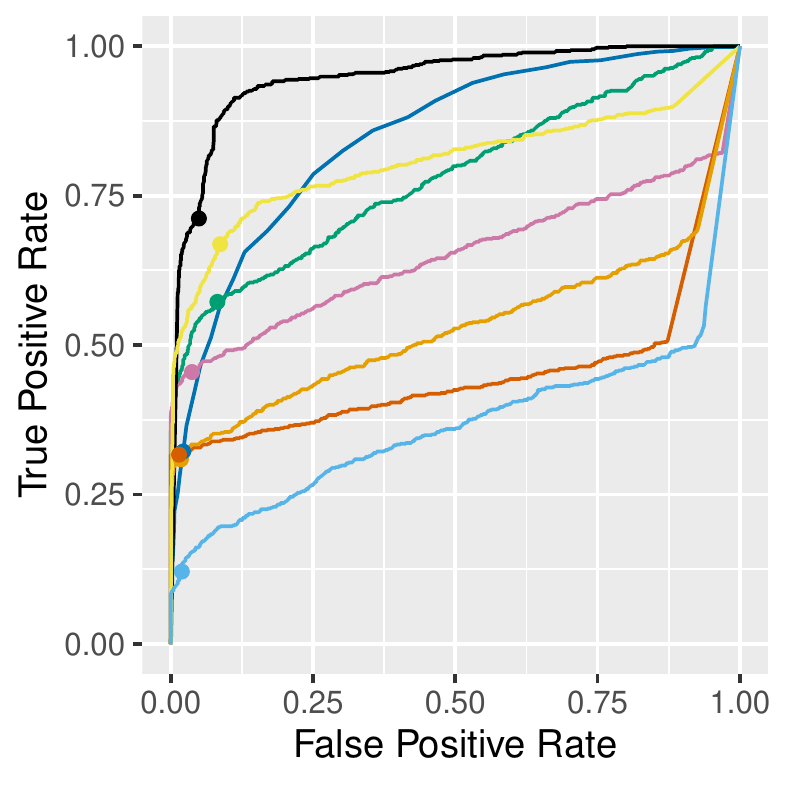}
        \vspace{-5mm}
        \caption{LCD triple $(C, X, Y)$}
        \label{fig:roc_lcd}
        \vspace{2mm}
    \end{subfigure}

    \begin{subfigure}[t]{0.245\textwidth}
        \centering
        \includegraphics[width=1\textwidth]{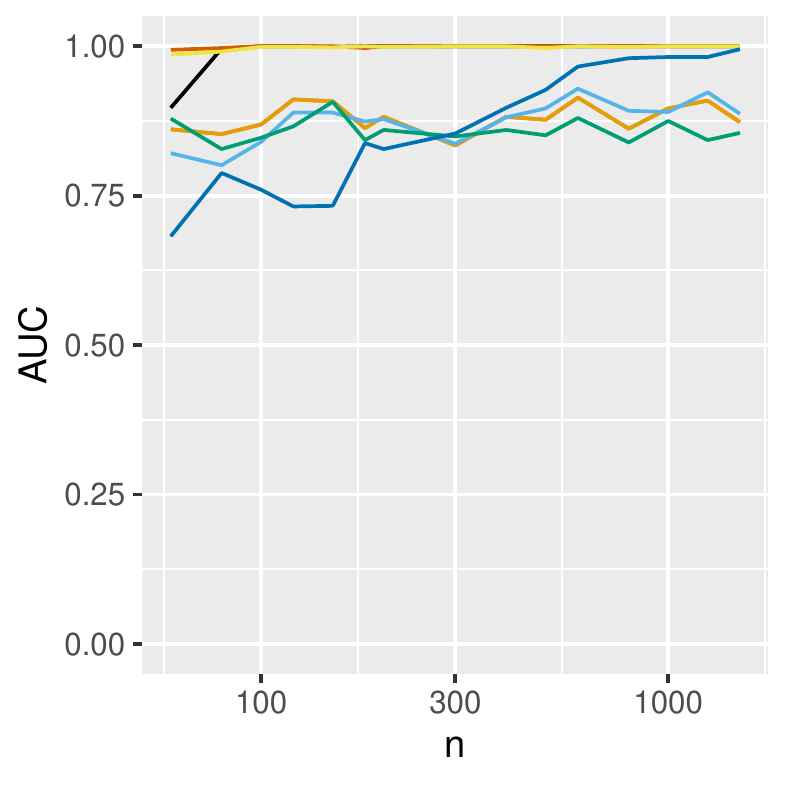}
        \vspace{-6mm}
        \caption{$C \dep X$}
        \label{fig:auc_ts}
    \end{subfigure}
    \begin{subfigure}[t]{0.245\textwidth}
        \centering
        \includegraphics[width=1\textwidth]{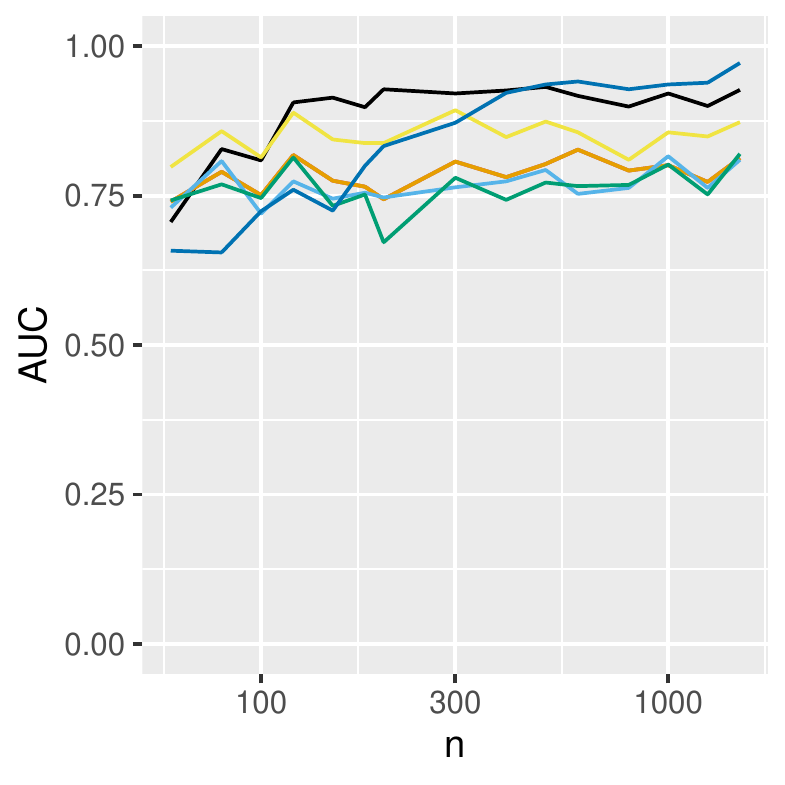}
        \vspace{-6mm}
        \caption{$X \dep Y$}
        \label{fig:auc_uci}
    \end{subfigure}
    \begin{subfigure}[t]{0.245\textwidth}
        \centering
        \includegraphics[width=1\textwidth]{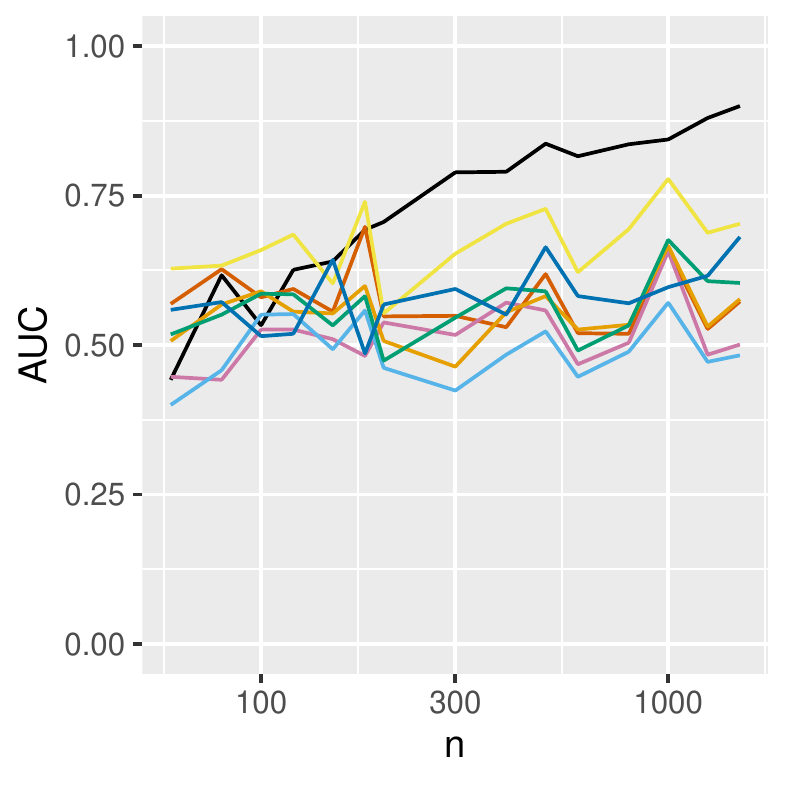}
        \vspace{-6mm}
        \caption{$C \dep Y|X$}
        \label{fig:auc_ci}
    \end{subfigure}
    \begin{subfigure}[t]{0.245\textwidth}
        \centering
        \includegraphics[width=1\textwidth]{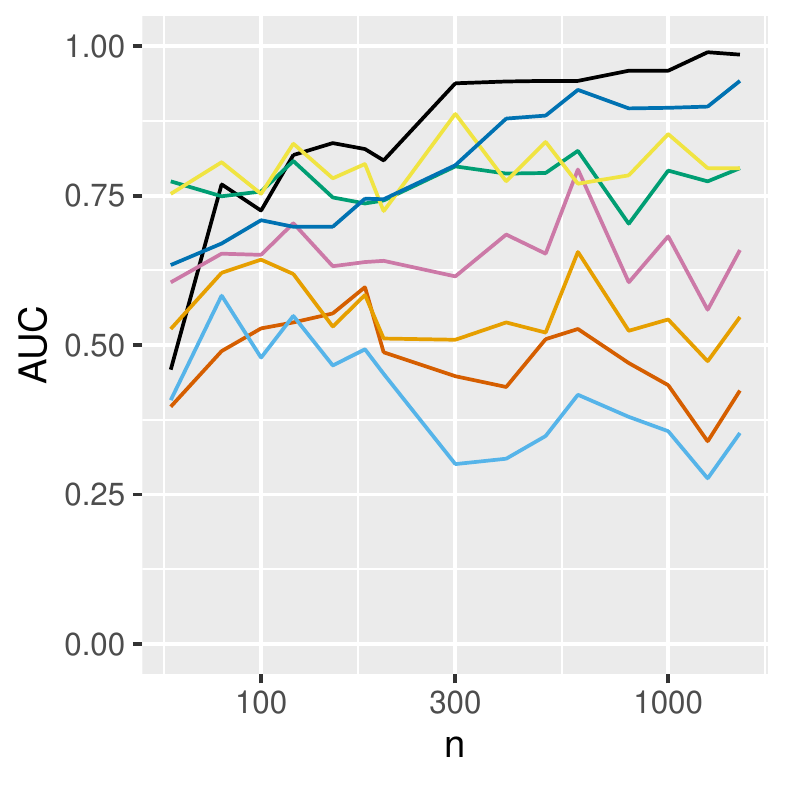}
        \vspace{-6mm}
        \caption{LCD triple $(C, X, Y)$}
        \label{fig:auc_lcd}
    \end{subfigure}

    \begin{subfigure}[t]{0.8\textwidth}
        \vspace{1mm}
        \centering
        \includegraphics[width=1\textwidth]{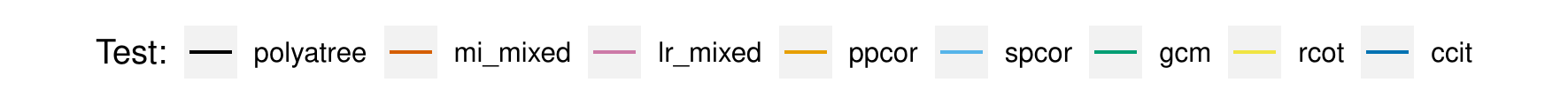}
    \end{subfigure}
    \vspace*{-2mm}
    \caption{ROC and AUC results for simulated data. The first row depicts ROC curves for individual tests (a--c) and for the LCD algorithm (d) over 2000 rounds of simulations at sample size $n=400$. The second row depicts the median AUC for varying sample size (ranging from 60 to 1500) for individual tests (e--g) and for the LCD algorithm (h)  over 200 rounds of simulations.}
    \label{fig:roc}
\end{figure*}

Figures \ref{fig:roc} (a--d) show the results of 2000 rounds of simulations, where in each round we simulate 400 observations. On the ROC curves we have marked the reference points $\alpha = 0.05$ and $\alpha=0.5$ for respectively frequentist and Bayesian tests. Figures \ref{fig:roc} (e--h) generalise these results, as they show the areas under the ROC curves (AUC) for varying sample sizes. We note that for conditional independence testing (Figure \ref{fig:roc_ci} and \ref{fig:auc_ci}), the P\'olya tree test from Section \ref{section:cond_two_sample_test} and the RCoT perform relatively well. It is interesting to see that the other tests have performance close to random guessing. It is however unclear whether this is due to the nonlinearity $\ell$, the intervention $g$ or the fact that $C$ is binary instead of continuous. From Figures \ref{fig:roc_lcd} and \ref{fig:auc_lcd} we see that the high performance of the P\'olya tree tests accumulates into good performance of the LCD algorithm. Interestingly, the CCIT also performs quite well, despite its weak performance in conditional independence testing.

\begin{figure}
    \centering
    \includegraphics[width=1\columnwidth]{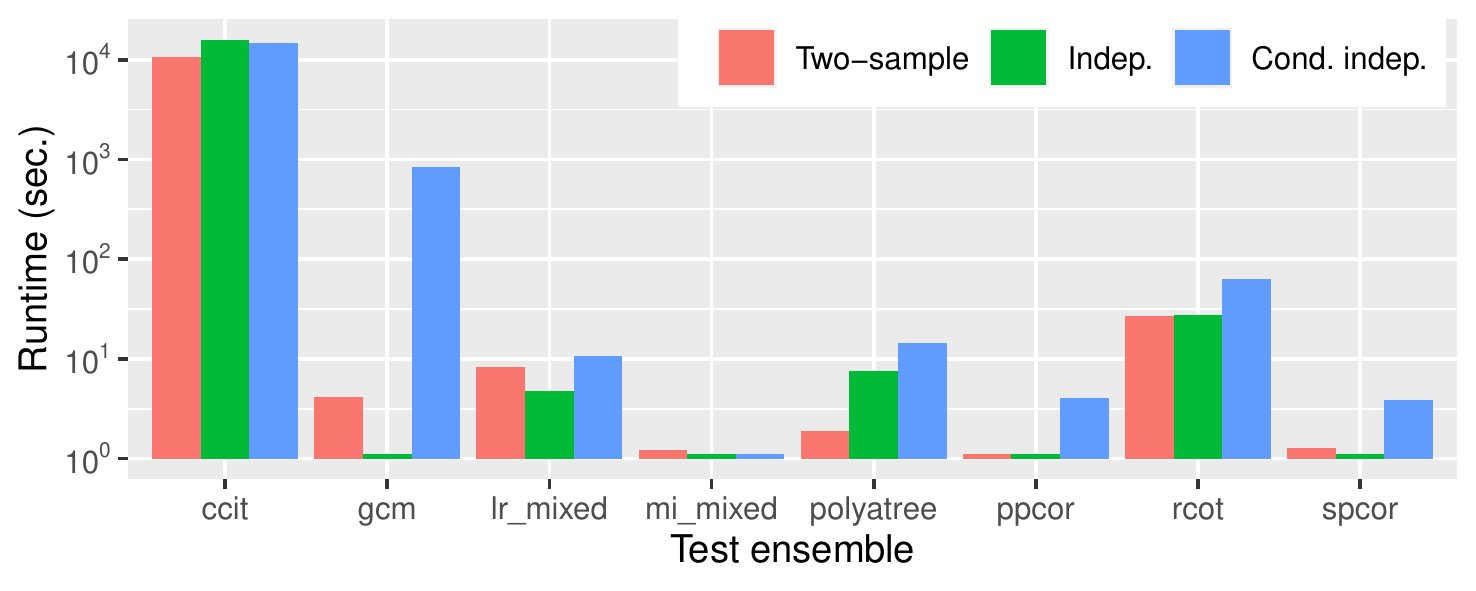}
    \caption{Runtimes of the different tests on 2000 rounds of simulations at $n=400$.}
    \label{runtimes}
\end{figure}

In Figure \ref{runtimes} we display for each independence test the computation times of the three test cases, accumulated over 2000 rounds of simulation at a sample size of $n=400$, as performed for generating Figures \ref{fig:roc} (a--c). The reader should be aware that for the GCM the difference in runtime between marginal and conditional independence testing is due to the fact that for conditional independence testing two nonlinear regressions are performed, and for marginal testing a statistic similar to partial correlation is computed. The CCIT has relatively high computation time due to costly training of the XGBoost classifier for each round of simulations, which makes it rather impractical to use. The partial correlation tests clearly perform best in terms of runtime. Overall, we conclude that the P\'olya tree tests provide a very good trade-off between performance and computation time.

\subsection{Protein expression data}
\label{section:sachs}
We apply the LCD algorithm, implemented with the Bayesian ensemble of independence tests, to protein expression data \citep{sachs2005}. \citeauthor{sachs2005} provide an `expert network', depicting the consensus (at that time) among biologists on the true network of signals between 11 proteins and phospholipids, and 10 reagents that are added to the cells. They estimate a causal graph which deviates from the expert network by some edges, refraining from claiming whether these edges should be added to the true network. For a detailed description of the data set and a depiction of the expert network we refer to the supplement.

Many authors have used this data set for estimating the underlying causal network, of which the graph of the original paper \citep{sachs2005} most closely resembles the expert network \citep{ramsey2018}. Furthermore, \citet{ramsey2018} and \citet{mooij2020} provide sufficient grounds for rejecting the expert network as being the true causal graph of the data. As we have no reliable ground truth to compare the output of the LCD algorithm with, we compare the output of LCD with its implementation with partial correlation.

\begin{figure*}
    \centering
    \begin{subfigure}{0.45\textwidth}
        \centering
        \includegraphics[width=1\textwidth]{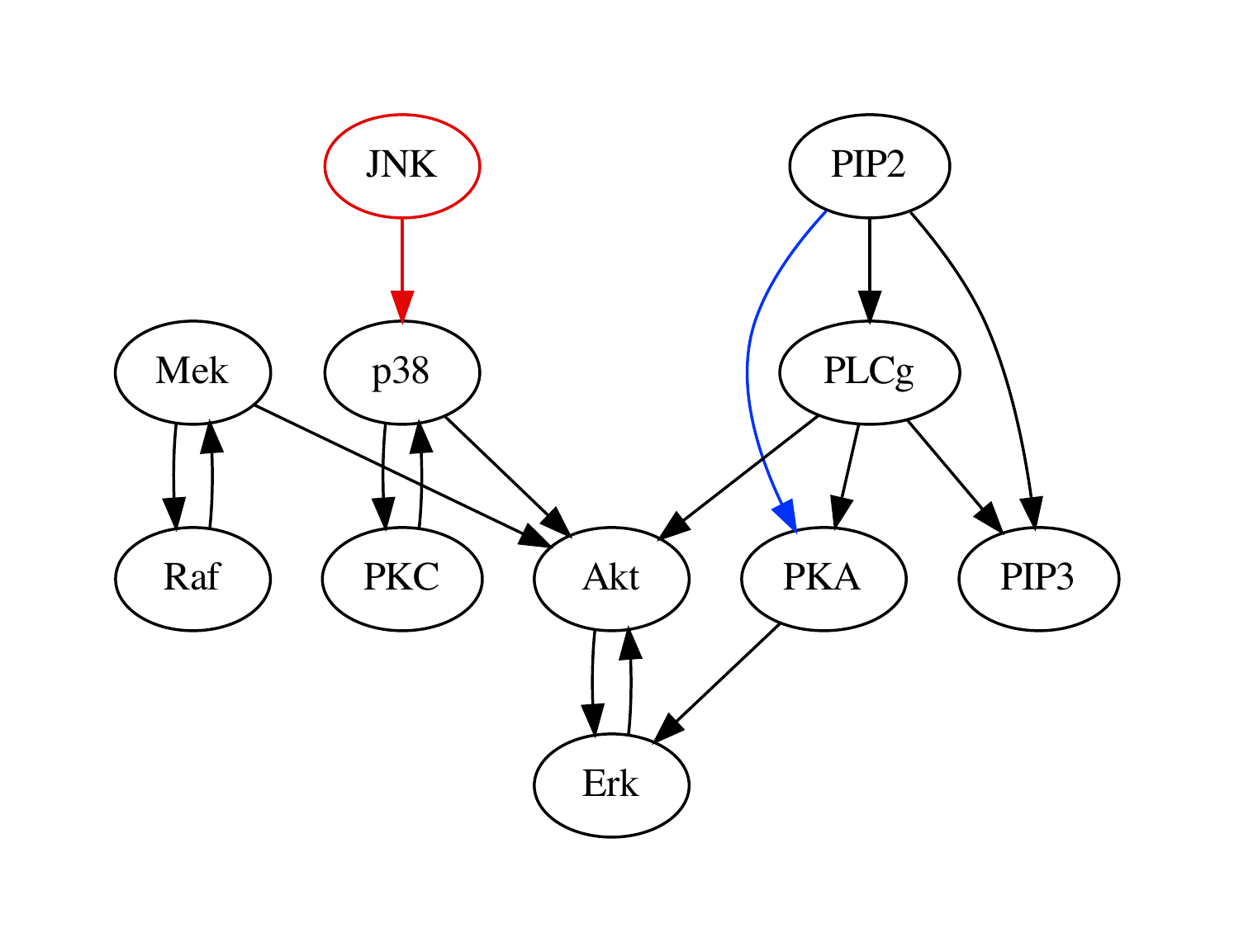}
        \caption{Output of LCD with P\'olya tree tests.}
        \label{fig:sachs_pt}
        \vspace*{2mm}
    \end{subfigure}
    \begin{subfigure}{0.45\textwidth}
        \centering
        \includegraphics[width=1\textwidth]{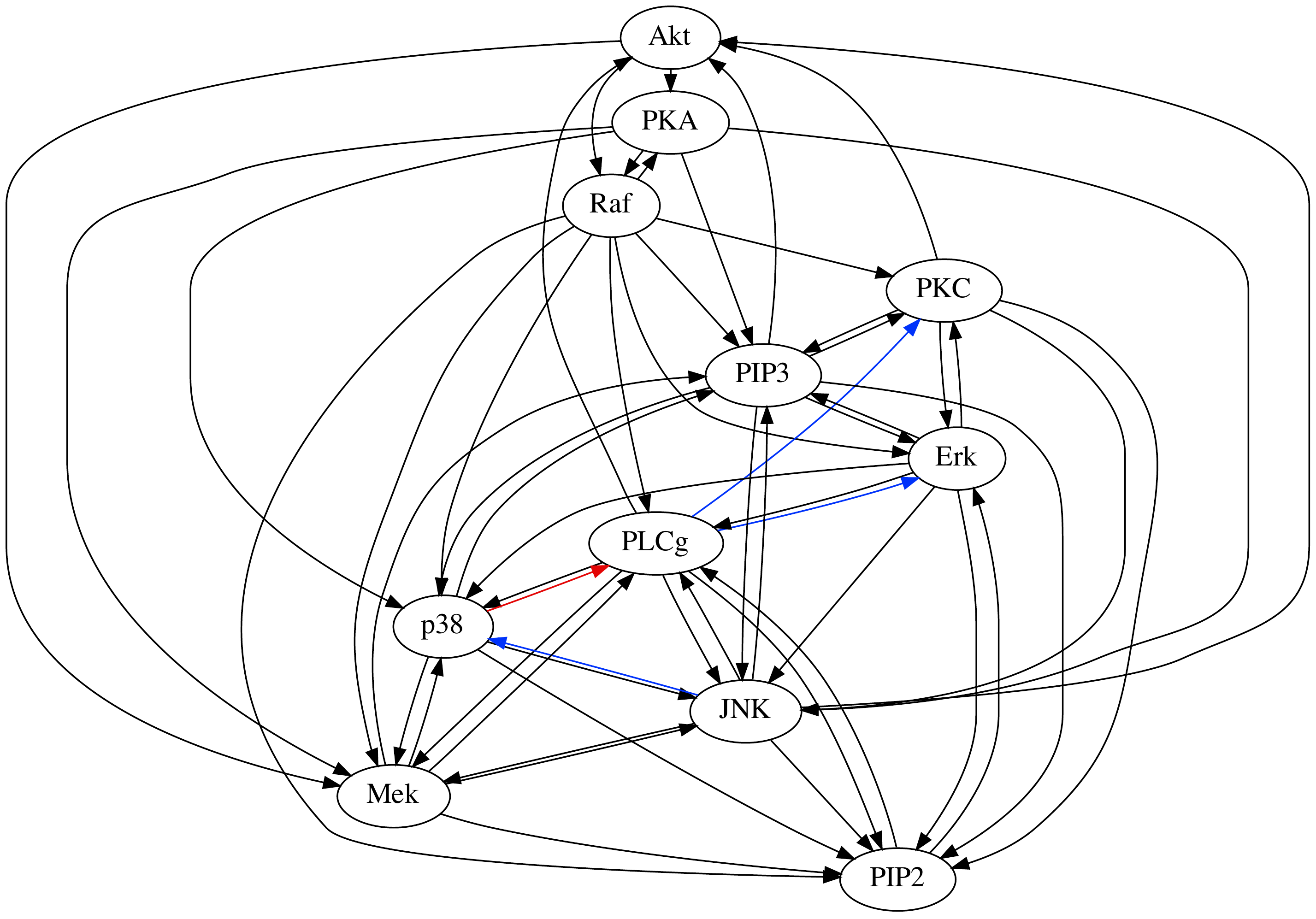}
        \caption{Output of LCD with the partial correlation test.}
        \label{fig:sachs_pcor}
    \end{subfigure}
    \caption{The output of LCD on the Sachs data. Edges indicate (possibly indirect) causal effects between the nodes. Black edges indicate strong evidence, red edges indicate substantial evidence, and blue edges indicate weak evidence.}
    \label{fig:sachs}
\end{figure*}

The output of the LCD algorithm implemented with the Bayesian tests and with the partial correlation test is shown in Figure \ref{fig:sachs}. In both cases we report the output of the LCD algorithm for multiple thresholds for the statistical tests. For the Bayesian tests (Figure \ref{fig:sachs_pt}) we use Bayes factor thresholds of $k=10$ (strong evidence, depicted in black), $k=4$ (substantial evidence, depicted in red) and $k=1$ (weak evidence, depicted in blue) \citep{kass1995}, and for the partial correlation test (Figure \ref{fig:sachs_pcor}) we report results for the p-value thresholds $\alpha = 0.0001$ (strong evidence, depicted in black), $\alpha = 0.005$ (substantial evidence, depicted in red) and $\alpha = 0.05$ (weak evidence, depicted in blue).

In general, we note that the output of LCD differs strongly among the use of different statistical tests, corroborating the premise that the performance of the algorithm highly depends on the choice of statistical test. Since the partial correlation test does not detect nonlinear conditional independencies, it has relatively low recall when compared with the P\'olya tree test, as shown in Figure \ref{fig:roc_ci}. This causes the LCD algorithm with partial correlations to output more false positives, resulting in a very dense causal graph, whereas LCD with the P\'olya tree tests produces a graph which is more likely to resemble the true causal model.

%% file: graphics/graph1.tex
\setlength{\fboxsep}{.8em}
\fbox{\begin{tikzpicture}[inner sep=1.5pt, -latex, on grid, auto, node distance=1.4cm, state/.style={circle, draw, minimum width=0.75cm}]
        \node[state] (C) {$C$};
        \node[state] (X) [right=of C] {$X$};
        \node[state] (Y) [below=of X] {$Y$};
        \node[state] (e_X) [right=of X] {$E_X$};
        \node[state] (e_Y) [right=of Y] {$E_Y$};
        \path (e_X) edge [] node[] {} (X);
        \path (C) edge [] node[] {} (X);
        \path (X) edge [] node[pos=0.4, left=0.3em] {$\ell$} (Y);
        \path (e_Y) edge [] node[] {} (Y);

        \node[] () at (C.west) {};
        \node[] () at (C.north) {};
        \node[] () at (e_X.east) {};

        \node[] () [below=2.3cm of Y, minimum height=3cm, label=center: {\footnotesize
        $\begin{aligned}
                X   & = g(C, E_X)                                    \\
                Y   & = \ell(X) + E_Y                                \\
                C   & \sim \textrm{Bernoulli}(1/2)                   \\
                E_X & \sim \mathcal{N}(0, 1)                         \\
                E_Y & \sim \mathcal{N}(0, \textrm{Var}(\ell(X)) / 4) \\
                    &                                                
            \end{aligned}$
        }] {};
    \end{tikzpicture}}

%% file: graphics/graph2.tex
\setlength{\fboxsep}{.8em}
\fbox{\begin{tikzpicture}[inner sep=1.5pt, -latex, on grid, auto, node distance=1.4cm, state/.style={circle, draw, minimum width=0.75cm}]
        \node[state] (C) {$C$};
        \node[state] (X) [right=of C] {$X$};
        \node[state] (Y) [below=of X] {$Y$};
        \node[state] (e_X) [right=of X] {$E_X$};
        \node[state] (e_Y) [right=of Y] {$E_Y$};
        \path (e_X) edge [] node[] {} (X);
        \path (C) edge [] node[] {} (X);
        \path (Y) edge [] node[pos=0.4, left=0.3em] {$\ell$} (X);
        \path (e_Y) edge [] node[] {} (Y);

        \node[] () at (C.west) {};
        \node[] () at (C.north) {};
        \node[] () at (e_X.east) {};

        \node[] () [below=2.3cm of Y, minimum height=3cm, label=center: {\footnotesize
        $\begin{aligned}
                X   & = g(C, \ell(Y) + E_X)                          \\
                Y   & = E_Y                                          \\
                C   & \sim \textrm{Bernoulli}(1/2)                   \\
                E_X & \sim \mathcal{N}(0, \textrm{Var}(\ell(Y)) / 4) \\
                E_Y & \sim \mathcal{N}(0, 1)                         \\
                    &
            \end{aligned}$
        }] {};
    \end{tikzpicture}}

%% file: graphics/graph3.tex
\setlength{\fboxsep}{.8em}
\fbox{\begin{tikzpicture}[inner sep=1.5pt, -latex, on grid, auto, node distance=1.4cm, state/.style={circle, draw, minimum width=0.75cm}]
        \node[state] (C) {$C$};
        \node[state] (X) [right=of C] {$X$};
        \node[state] (Y) [below=of X] {$Y$};
        \node[state] (L) [left=of Y] {$L$};
        \node[state] (e_X) [right=of X] {$E_X$};
        \node[state] (e_Y) [right=of Y] {$E_Y$};
        \path (e_X) edge [] node[] {} (X);
        \path (C) edge [] node[] {} (X);
        \path (L) edge [] node[pos=0.5] {$\ell$} (X);
        \path (L) edge [] node[pos=0.4, above=0.2em] {$\ell$} (Y);
        \path (e_Y) edge [] node[] {} (Y);

        \node[] () at (C.west) {};
        \node[] () at (C.north) {};
        \node[] () at (e_X.east) {};

        \node[] () [below=2.3cm of Y, minimum height=3cm, label=center: {\footnotesize
        $\begin{aligned}
                X   & = g(C, \ell(L) + E_X)                        \\
                Y   & = \ell(L) + E_Y                              \\
                C   & \sim \textrm{Bernoulli}(1/2)                 \\
                L   & \sim \mathcal{N}(0, 1)                       \\
                E_X & \sim \mathcal{N}(0, \textrm{Var}(\ell(L))/4) \\
                E_Y & \sim \mathcal{N}(0, \textrm{Var}(\ell(L))/4)
            \end{aligned}$
        }] {};
    \end{tikzpicture}}

%% file: input/4_sensitivity_analysis.tex
\begin{figure*}[h]
    \centering
    \begin{subfigure}[t]{0.245\textwidth}
        \centering
        \includegraphics[width=1\textwidth]{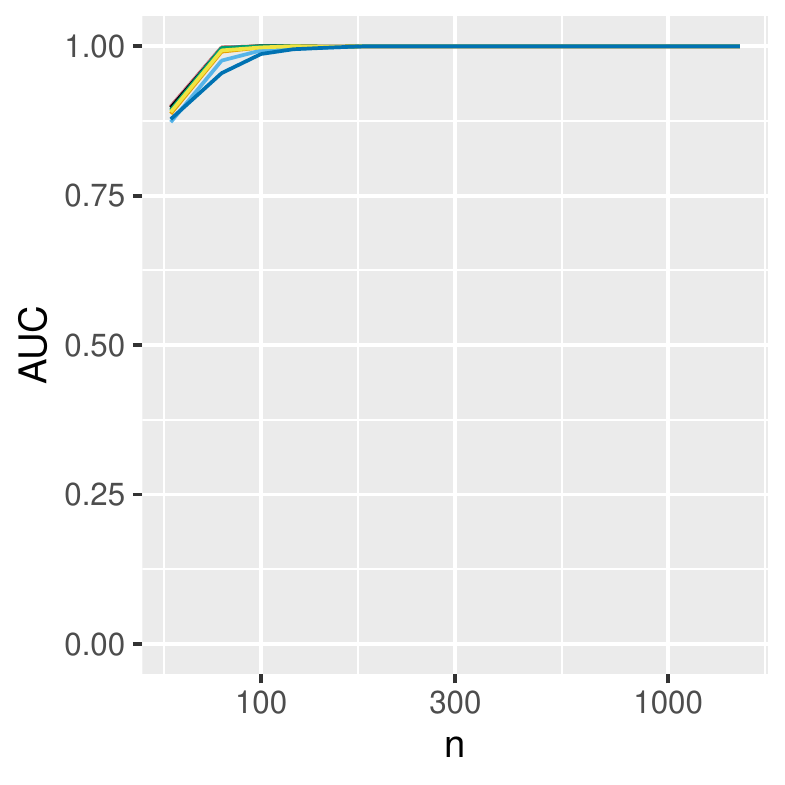}
        \vspace{-6mm}
        \caption{$C \dep X$}
    \end{subfigure}
    \begin{subfigure}[t]{0.245\textwidth}
        \centering
        \includegraphics[width=1\textwidth]{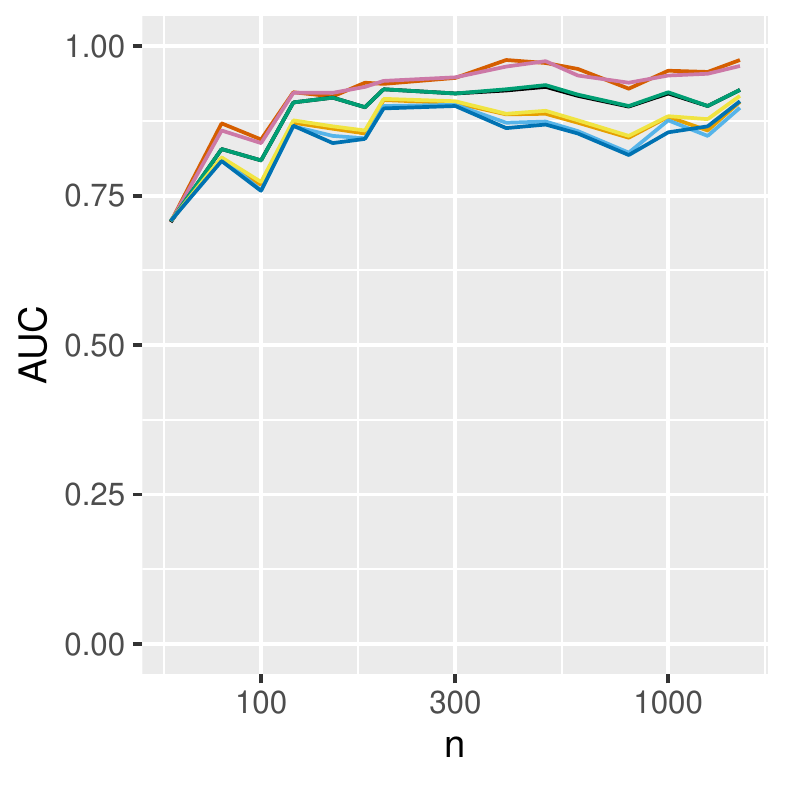}
        \vspace{-6mm}
        \caption{$X \dep Y$}
    \end{subfigure}
    \begin{subfigure}[t]{0.245\textwidth}
        \centering
        \includegraphics[width=1\textwidth]{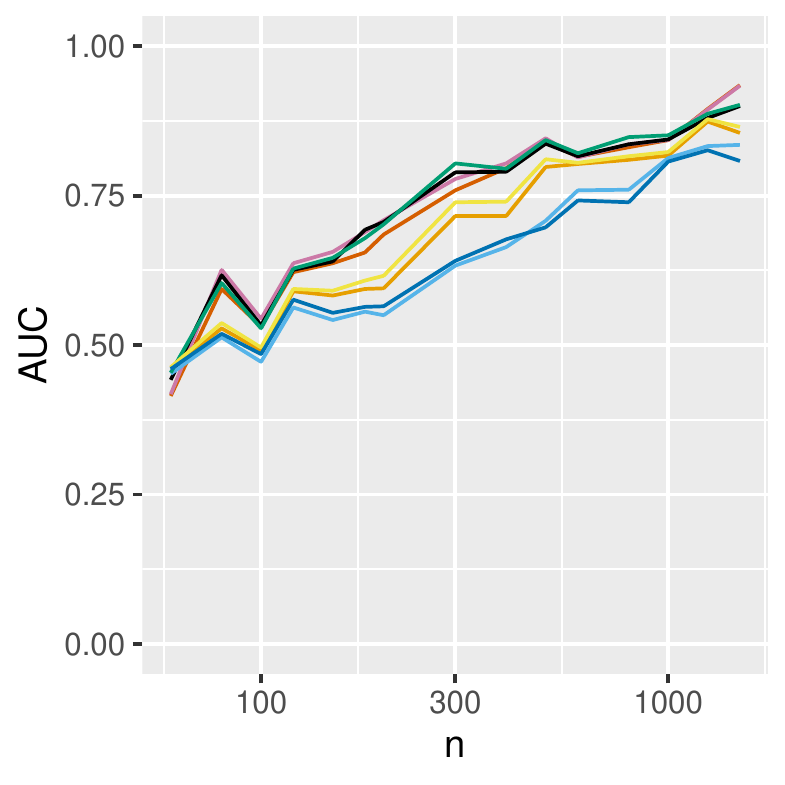}
        \vspace{-6mm}
        \caption{$C \dep Y|X$}
    \end{subfigure}
    \begin{subfigure}[t]{0.245\textwidth}
        \centering
        \includegraphics[width=1\textwidth]{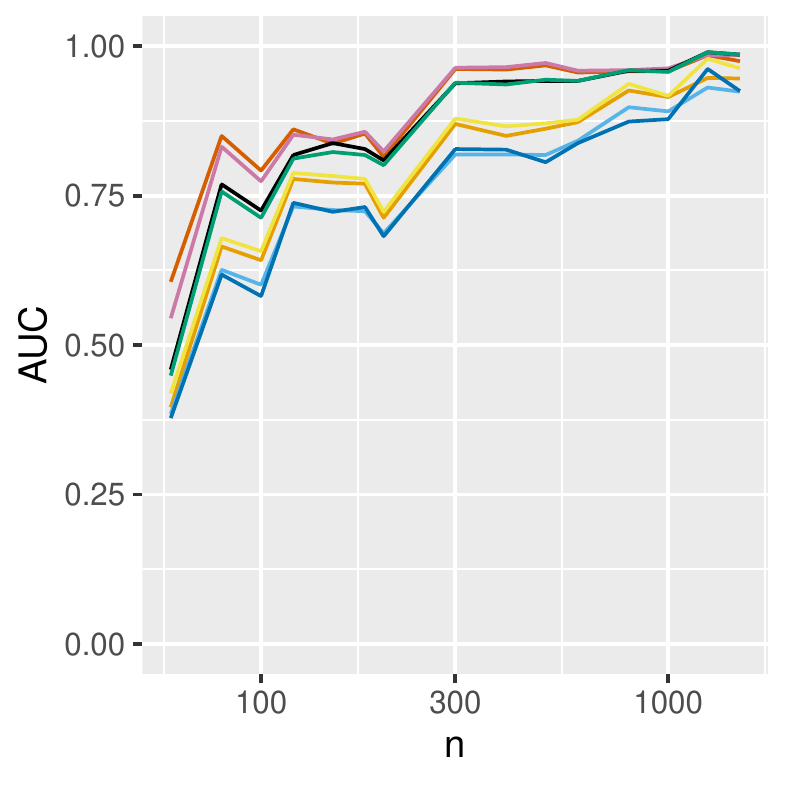}
        \vspace{-6mm}
        \caption{LCD triple $(C, X, Y)$}
    \end{subfigure}

    \vspace{3mm}
    \begin{subfigure}[t]{0.8\textwidth}
        \centering
        $\rho(j):$\hspace{1pt}
        \vcenteredhbox{\includegraphics[width=4.6mm]{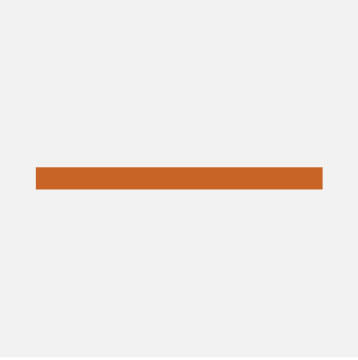}}
        \hspace{1pt}$\tfrac{1}{10}j^2$\hspace{1pt}
        \vcenteredhbox{\includegraphics[width=4.6mm]{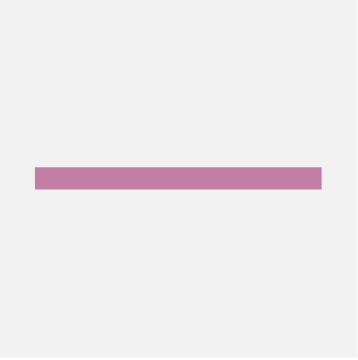}}
        \hspace{1pt}$\tfrac{1}{5}j^2$\hspace{1pt}
        \vcenteredhbox{\includegraphics[width=4.6mm]{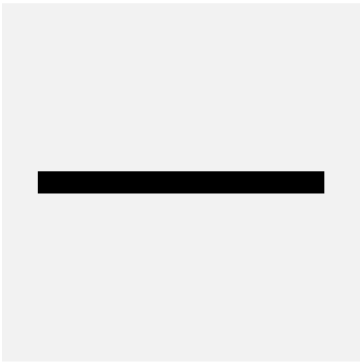}}
        \hspace{1pt}$j^2$\hspace{1pt}
        \vcenteredhbox{\includegraphics[width=4.6mm]{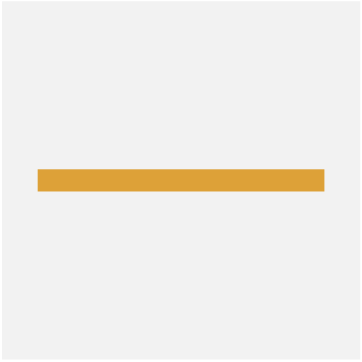}}
        \hspace{1pt}$5j^2$\hspace{1pt}
        \vcenteredhbox{\includegraphics[width=4.6mm]{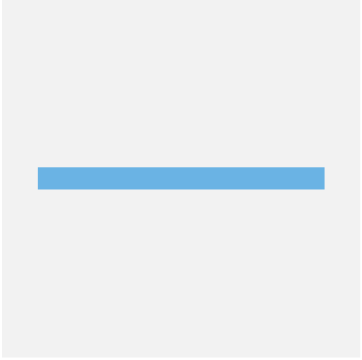}}
        \hspace{1pt}$10j^2$\hspace{1pt}
        \vcenteredhbox{\includegraphics[width=4.6mm]{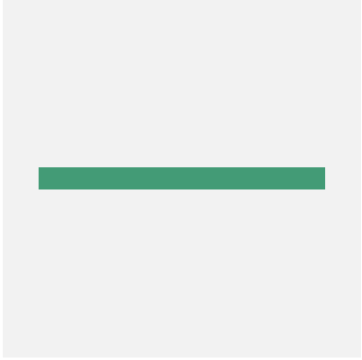}}
        \hspace{1pt}$2^{j}$\hspace{1pt}
        \vcenteredhbox{\includegraphics[width=4.6mm]{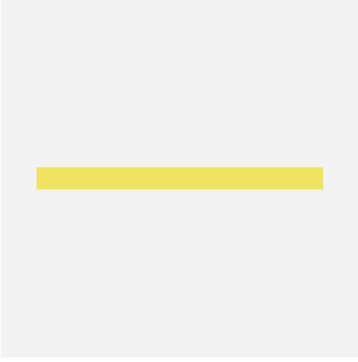}}
        \hspace{1pt}$4^{j}$\hspace{1pt}
        \vcenteredhbox{\includegraphics[width=4.6mm]{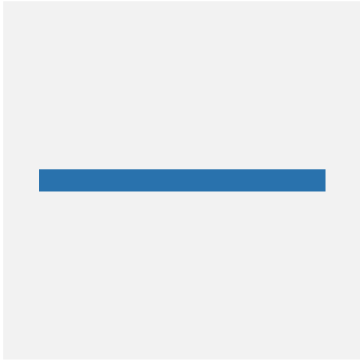}}
        \hspace{1pt}$8^{j}$\hspace{1pt}
    \end{subfigure}
    \caption{Sensitivity of the performance of the P\'olya tree tests with respect to the parameter $\mathcal{A}$.}
    \label{fig:sensitivity}
\end{figure*}

As mentioned earlier, the P\'olya tree is parametrised by the set $\mathcal{A}$, where in the previous section we have used $\alpha_\kappa = |\kappa|^2$. In general we can let $\alpha_\kappa := \rho(|\kappa|)$ for any positive function $\rho$, in which case we have
\begin{equation}
    \textrm{Var}(\mathcal{P}(B_\kappa)) = \frac{1}{4^{|\kappa|}}\left(\prod_{j=1}^{|\kappa|}\frac{2\rho(j) + 2}{2\rho(j) + 1} - 1\right),
\end{equation}
and samples from the P\'olya tree are dominated by Lebesgue measure if $\sum_{j=1}^\infty \rho(j)^{-1} < \infty$ \citep{kraft1964}. \citeauthor{walker1999} propose to use $\rho(j) = cj^2$ for some $c>0$, in which case decreasing $c$ increases the variance of $\mathcal{P}$, causing $\mathcal{P}$ to be less dependent on the choice of $G$. We have chosen $c=1$ as a default value in the previous section as it is promoted as a “sensible canonical choice” by \citet{lavine1992}. According to \citet{holmes2015}, having $c$ between 1 and 10 is in general a good choice. To obtain an better understanding of the dependency of the P\'olya tree on this parameter, we have repeated the experiments of Figure \ref{fig:roc} (e--h) for different choices of $\rho$. More specifically, we have repeated the experiments for parameters $\rho(j) = \tfrac{1}{10}j^2, \tfrac{1}{5}j^2, j^2, 5j^2, 10j^2, 2^{j}, 4^{j}$ and $8^{j}$ \citep{berger2001}. The results are shown in Figure \ref{fig:sensitivity}. We note that the performance of the tests is not heavily influenced by the choice of $\mathcal{A}$, and that $\rho(j) = j^2$ seems to be an appropriate default value.

%% file: input/5_discussion.tex
In this work we have proposed a novel nonparametric conditional two-sample test, which is possibly the first conditional independence test of this type. The test is analysed in its own right and as a subroutine of the Local Causal Discovery algorithm, and in both cases can outperform current state-of-the-art nonparametric continuous conditional independence tests and parametric mixed conditional independence tests. However, we have made some modelling decisions which might be reconsidered when using this test in practice.

First we note that the choice of $\mathcal{A}$ may influence the suitability of the test. Section \ref{section:sensitivity} suggests that $\alpha_\kappa = |\kappa|^2$ is a sensible parametrisation, but this may be reconsidered in applications. Another consideration is the choice of the family of partitions $\mathcal{T}$. Having canonical partitions increases the intelligibility of the P\'olya tree, but essentially any recursive partitioning tree suffices. We note that the maximum partitioning depth $J = \lfloor\log_4(n)\rfloor$ is quite arbitrarily chosen to reduce computation time. However, as our choice of $\alpha_\kappa$ implies relatively low dependence on the base measure $G$ and as we standardise the data to approximately fit the standard Gaussian base measure, we believe that we have chosen sensible default parameters.

In general, it is hard to theoretically analyse for which types of distributions conditional independence tests work properly. For frequentist tests, the asymptotic distribution of the test statistic is often provided, which holds under rather technical assumptions which may be hard to validate against a provided dataset (see \cite{strobl2019} for an example of such assumptions). The same holds for theoretical consistency results of the test statistic under the alternative. \cite{shah2020} show that in order to have power against an acceptably large set of alternatives, one should restrict the set of distributions considered under $H_0$. In a Bayesian setting, consistency of the Bayes Factor is determined by whether the posterior corresponding to the true hypothesis is consistent (i.e. the marginal likelihood is large), and the marginal likelihood remains small under the false hypothesis. Sufficient conditions for posterior convergence are e.g. provided by Doob’s Theorem and Schwartz’s Theorem \citep{ghosal2017}, but necessary conditions (which could be used to restrict $H_0$ and $H_1$) are not available to our knowledge. One should also investigate the behaviour of the posterior likelihood under misspecification to properly determine for which $H_0$ and $H_1$ the test works properly.

Many constraint-based causal inference algorithms (other than LCD) require conditional independence testing of the form $C \indep X|Z$ for $d$-dimensional $Z$ with $d>1$. Extending our method is straightforward, as the canonical partitions of $\mathcal{Z}$ can be constructed as the per-level cartesian product of $d$ one-dimensional canonical partitions \citep{hanson2006}. However, this extension suffers from the curse of dimensionality, so further research should look into how this problem can be mitigated.

This work only addresses testing $C\indep X|Z$ where $C$ is binary. Although this test is already of high importance to the field of causal discovery, extending this test to discrete $C$ would be of real use and is the subject of current research.

The ensemble of P\'olya tree prior based independence tests provides good results when utilised in a causal inference algorithm applied on synthetic data, and produces sensible output on real world data. We therefore believe that it is a promising area of research, which hopefully will improve the robustness and applicability of causal inference algorithms.

%% file: input/A_polyatrees.tex
In general, our setup for independence testing will assume availability of independent samples $X_1, ..., X_n$ of a random variable $X$ with continuous distribution $P$. We let $\mathcal{X}$ denote the domain of $X$, and let $\mathcal{M}$ be the space of continuous distributions on $\mathcal{X}$. Our hypotheses will be of the form
\begin{equation}\label{hypotheses}
    H_0: X \sim P \textrm{ with } P\in \mathcal{M}_0, \quad H_1: X \sim P \textrm{ with } P\in \mathcal{M}_1,
\end{equation}
where $\mathcal{M}_0, \mathcal{M}_1 \subset \mathcal{M}$, and $\mathcal{M}_0 \cap \mathcal{M}_1 = \emptyset$. Since we wish to device a Bayesian test, we will define prior distributions $\Pi_0$ and $\Pi_1$ with support on $\mathcal{M}_0$ and $\mathcal{M}_1$ respectively. Then we compare the evidence of the models given the data via the Bayes factor, i.e.
\begin{equation} \label{BF_def}
    \mathrm{BF}_{01} = \frac{\PP(H_0 | X_{1:n})}{\PP(H_1 | X_{1:n})} = \frac{p(X_{1:n} | H_0)}{p(X_{1:n}|H_1)}\frac{\PP(H_0)}{\PP(H_1)} = \frac{\int_{\mathcal{M}}\prod_{i=1}^n p(X_i)d\Pi_0(P)}{\int_{\mathcal{M}}\prod_{i=1}^n p(X_i)d\Pi_1(P)}
\end{equation}
where we have placed equal prior weights on $H_0$ and $H_1$, so $\PP(H_0) = \PP(H_1) = 1/2$.

A canonical choice for a prior on a space of probability distributions is the Dirichlet Process. However, samples from the Dirichlet process are almost surely discrete distributions, so the Dirichlet Process is not a suitable choice for our setup. The P\'olya tree prior does not suffer from this characteristic \citep{ferguson1974}, and can be parametrised to be a suitable prior on $\mathcal{M}$. Since the elements of $\mathcal{M}$ have support on $\mathcal{X}$, we will speak of a P\'olya tree on $\mathcal{X}$. We will first construct a P\'olya tree on $\mathcal{X} \subseteq \RR$, and then extend this definition to a P\'olya tree on $\mathcal{X}\times\mathcal{Y} \subseteq \RR^2$.

First we recall the construction of the one-dimensional P\'olya tree as described in the main paper. In particular, we construct a P\'olya tree on $(\mathcal{X}, \mathcal{B}(\mathcal{X}))$, where $\mathcal{X}\subseteq \RR$, and $\mathcal{B}(\mathcal{X})$ denotes the Borel sigma-algebra on $\mathcal{X}$. In order to construct a random measure on $\mathcal{B}(\mathcal{X})$, we will assign random probabilities to a family of subsets $\mathcal{T}$ of $\mathcal{X}$ which generates the Borel sets. The family of subsets that we consider are the dyadic partitions of $[0,1]$, mapped under the inverse of some cumulative distribution function $G$ on $\mathcal{X}$. This results in the \textit{canonical} family of partitions of $\mathcal{X}$, where for level $j$ we have $\mathcal{X} = \bigcup_{\kappa \in \{0,1\}^j} B_\kappa$, with
\begin{equation}\label{app:1d_partition}
    B_\kappa := [G^{-1}(\tfrac{k-1}{2^j}), G^{-1}(\tfrac{k}{2^j})),
\end{equation}
and $k$ is the natural number corresponding to the bit string $\kappa \in \{0,1\}^j$. A schematic depiction of this binary tree of partitions is shown in Figure \ref{app:1d_partition_diagram}. We define the index set by $K := \{\{0,1\}^j: j\in \NN\}$, so the family of subsets of $\mathcal{X}$ that we consider is $\mathcal{T} := \{B_\kappa : \kappa\in K\}$. From basic measure theory we know that $\mathcal{T}$ indeed generates $\mathcal{B}(\mathcal{X})$. We assign random probabilities to the elements of $\mathcal{T}$ by first assigning random probabilities to $B_0$ and $B_1$, and randomly subdividing these masses among the children of $B_0$ and $B_1$. In particular, for the first level of the partition we assign the random probabilities $\mathcal{P}(B_0) = \theta_0$ and $\mathcal{P}(B_1) = \theta_1$ with $(\theta_0, \theta_1) \sim \textrm{Dir}(\alpha_0, \alpha_1)$, for some hyper-parameters $\alpha_0$ and $\alpha_1$. Then, for every $B_\kappa\in \mathcal{T}$ we split the mass that is assigned to $B_\kappa$ by assigning a fraction $\theta_{\kappa 0}$ to $B_{\kappa 0}$ and a fraction $\theta_{\kappa 1}$ to $B_{\kappa 1}$, where we let $(\theta_{\kappa 0}, \theta_{\kappa 1}) \sim \textrm{Dir}(\alpha_{\kappa 0}, \alpha_{\kappa 1})$. This construction yields a P\'olya tree on $\mathcal{X}$, which is a random measure on $\mathcal{B}(\mathcal{X})$:

\begin{figure}[htb]
    \centering
    \fbox{\input{graphics/1d_partition}}
    \caption{Construction of a one-dimensional P\'olya tree based on canonical partitions.}
    \label{app:1d_partition_diagram}
\end{figure}

\begin{definition}[\citealp{lavine1992}]
A random probability measure $\mathcal{P}$ on $(\mathcal{X}, \mathcal{B}(\mathcal{X}))$ is said to have a P\'olya tree distribution with parameter $(\mathcal{T}, \mathcal{A})$, written $\mathcal{P} \sim {\normalfont\textrm{PT}}(\mathcal{T}, \mathcal{A})$, if there exist nonnegative numbers $\mathcal{A} = \{(\alpha_{\kappa 0},\alpha_{\kappa 1}) : \kappa \in K\}$ and random variables $\Theta = \{(\theta_{\kappa 0}, \theta_{\kappa 1}) : \kappa \in K\}$ such that the following hold:
\begin{enumerate}[itemsep=1.8pt, parsep=1.8pt]
    \item all the random variables in $\Theta$ are independent;
    \item for every $\kappa \in K$, we have $(\theta_{\kappa 0}, \theta_{\kappa 1}) \sim {\normalfont\textrm{Dir}}(\alpha_{\kappa 0}, \alpha_{\kappa 1})$;
    \item for every $j \in \NN$ and every $\kappa \in \{0,1\}^j$ we have $\mathcal{P}(B_{\kappa} | \Theta) = \prod_{i=1}^j \theta_{\kappa_1 ...\kappa_{i}}$.
\end{enumerate}
\end{definition}

The support of the P\'olya tree is determined by the choice of $\mathcal{T}$ and $\mathcal{A}$. In general, any separating binary tree of partitions of $\mathcal{X}$ can be considered. In this paper we only consider partitions of the type of equation (\ref{app:1d_partition}). \cite{ferguson1974} shows that the P\'olya tree is a Dirichlet process if $\alpha_\kappa = \alpha_{\kappa 0} + \alpha_{\kappa 1}$. The parameter of this Dirichlet process is the mean of the P\'olya tree, i.e.\ the probability measure $G_0$ on $\mathcal{B}(\mathcal{X})$ defined by $G_0(B) := \EE(\mathcal{P}(B))$ for all $B\in\mathcal{B}(\mathcal{X})$ \citep{lavine1994}. This implies that for this choice of $\mathcal{A}$, the support of the P\'olya tree is contained in the space of discrete distributions. Sufficient conditions on $\mathcal{A}$ for samples of the P\'olya tree to be continuous distributions are given by the following theorem:

\begin{theorem}[\citealp{kraft1964}]\label{thm_kraft} Let $\bar{\sigma}_j := \sup\{\textrm{Var}(\theta_\kappa) : \kappa \in \{0,1\}^j\}$. If $\mathbb{E}(\theta_\kappa) = 1/2$ for all $\kappa \in K$ and $\sum_{j=1}^\infty \bar{\sigma}_j < \infty$, then with probability one, samples from $\mathcal{P}$ are absolutely continuous with respect to Lebesgue measure.
\end{theorem}

This condition is satisfied if for each $\kappa\in \{0,1\}^j$ we take $\alpha_{\kappa 0} = \alpha_{\kappa 1} = j^2$, which is promoted as a ‘sensible canonical choice’ by \citet{lavine1992}. In this case we indeed have $\EE(\theta_\kappa) = 1/2$, and thus for every $j\in \NN$, the mass is (in expectation) split uniformly over the $B_\kappa$ for all $\kappa\in\{0,1\}^j$. As a consequence the P\'olya tree is centred on the base distribution with cumulative distribution function $G$, i.e.\ $\EE(\mathcal{P}(B_\kappa)) = \int_{B_\kappa} G'(x)dx$. As mentioned in the main paper we only consider partitions up to a pre-determined level $J(n)$.

Let $X$ be a continuous random variable with a distribution that lies in the support of the P\'olya tree $\mathcal{P} \sim \textrm{PT}(\mathcal{T}, \mathcal{A})$. Drawing a distribution from $\mathcal{P}$ is done by drawing from each of the random variables in $\Theta$. If we let $X_1, ..., X_n$ be a sample from $X$, then the likelihood of that sample with respect to a sampled distribution $\Theta$ from the P\'olya tree $\textrm{PT}(\mathcal{T}, \mathcal{A})$ is
\begin{equation}
    p (X_{1:n}|\Theta, \mathcal{T}, \mathcal{A}) = \prod_{\kappa \in K}\theta_{\kappa 0}^{n_{\kappa 0}}\theta_{\kappa 1}^{n_{\kappa 1}},
\end{equation}
where $n_\kappa$ denotes the number of observations lying in $B_\kappa$, i.e.\ $n_\kappa := |X_{1:n}\cap B_\kappa|$. If we integrate over all possible values of all $\theta_\kappa$, we obtain the marginal likelihood
\begin{equation}\label{app:1d_likelihood}
    p(X_{1:n} | \mathcal{T}, \mathcal{A}) = \prod_{\kappa \in K}\frac{\textrm{B}(\alpha_{\kappa 0} + n_{\kappa 0}, \alpha_{\kappa 1} + n_{\kappa 1})}{\textrm{B}(\alpha_{\kappa 0}, \alpha_{\kappa 1})},
\end{equation}
where $\textrm{B}(\cdot)$ denotes the Beta function. Note that this quantity corresponds to the marginal likelihood $\int_{\mathcal{M}}\prod_{i=1}^np(X_i)d\Pi(P)$, a version of which occurs in the numerator and denominator of the right-hand side of equation (\ref{BF_def}). This marginal likelihood will therefore be a fundamental quantity in the Bayesian tests that we consider.

\subsection{A nonparametric two-sample test}
\label{two_sample_test}
In order to use the P\'olya tree prior for Bayesian testing, we have to formulate our hypotheses $H_0$ and $H_1$ in terms of the relevant spaces of distributions $\mathcal{M}_0$ and $\mathcal{M}_1$, as suggested by equation (\ref{hypotheses}). This is done by picking P\'olya tree prior $\mathcal{P}_i\sim\Pi_i$ under $H_i$, and defining $\mathcal{M}_i$ to be the support of $\Pi_i$, for $i = 0, 1$. Given data to test our hypothesis with, we calculate marginal likelihoods via equation (\ref{app:1d_likelihood}) for both P\'olya trees $\mathcal{P}_0$ and $\mathcal{P}_1$, which are in turn used for calculating the Bayes factor via (\ref{BF_def}).

We first use this procedure to describe the nonparametric two-sample test, as proposed by \citet{holmes2015}. Given a sample $\{(X_1, B_1), ..., (X_n, B_n)\}$ from binary variable $C$ and continuous variable $X$, define $X^{(0)} := \{X_i : C_i = 0, i =1, .., n\}$ and $X^{(1)} := \{X_i : C_i = 1, i = 1, .., n\}$. Let $F$ denote the distribution of $X$, and let $F^{(0)}$ and $F^{(1)}$ denote the distributions of $X^{(0)}$ and $X^{(1)}$. We formulate the independence between $X$ and $C$ as a two-sample test, i.e.
\begin{align}
    & H_0: X \indep C \iff F^{(0)} = F^{(1)} = F  \\
    & H_1: X \dep C \iff F^{(0)} \neq F^{(1)}.
\end{align}

Under $H_0$ we standardise the sample $X_{1:n}$, and compute its marginal likelihood using equation (\ref{app:1d_likelihood}). Under $H_1$, we model $X^{(0)}$ and $X^{(1)}$ as being samples from independent random variables, having different distributions. Since separately normalising $X^{(0)}$ and $X^{(1)}$ may erase distinctive features between the samples, we first standardise $X$, and then subdivide $X$ into $X^{(0)}$ and $X^{(1)}$.

We formulate the Bayes factor as
\begin{equation}\label{two_sample_BF}
    \mathrm{BF}_{01} = \frac{p(X_{1:n} | \mathcal{T}, \mathcal{A})}{p(X^{(0)} | \mathcal{T}, \mathcal{A})p(X^{(1)} | \mathcal{T}, \mathcal{A})}.
\end{equation}

Upon inspection of equation (\ref{app:1d_likelihood}) we see that the Bayes factor can be written as an infinite product of fractions, being
\begin{align}
\begin{split}
    \mathrm{BF}_{01} = \prod_{\kappa\in K}
    & \frac{\textrm{B}(\alpha_{\kappa 0} + n_{X|\kappa 0}, \alpha_{\kappa 1} + n_{X|\kappa 1})\textrm{B}(\alpha_{\kappa 0}, \alpha_{\kappa 1})}{\textrm{B}(\alpha_{\kappa 0} + n_{X^{(0)}|\kappa 0}, \alpha_{\kappa 1} + n_{X^{(0)}|\kappa 1})\textrm{B}(\alpha_{\kappa 0} + n_{X^{(1)}|\kappa 0}, \alpha_{\kappa 1} + n_{X^{(1)}|\kappa 1})},
\end{split}
\end{align}
where $n_{X|\kappa} := |X_{1:n}\cap B_\kappa|$, and $n_{X^{(0)}|\kappa}$, $n_{X^{(1)}|\kappa}$ are defined similarly. We note that whenever $n_{X|\kappa} \leq 1$ the fraction has a value of 1, so we calculate the marginal likelihoods until we either reach the maximum partitioning depth $J(n)$, or until $n_{X|\kappa}\leq1$.

\subsection{Two-dimensional P\'olya trees}
Now that we have defined a P\'olya tree on $(\mathcal{X}, \mathcal{B}(\mathcal{X}))$ with $\mathcal{X} \subseteq \RR$, we extend this definition to a P\'olya tree on $(\mathcal{X}\times\mathcal{Y}, \mathcal{B}(\mathcal{X} \times \mathcal{Y}))$ with $\mathcal{X}\times\mathcal{Y} \subseteq \RR^2$. This construction is done similarly to the construction on $\mathcal{X}$. We consider a base measure with cumulative distribution function $G$ on $\mathcal{X} \cup \mathcal{Y}$, and partition $\mathcal{X} \times \mathcal{Y}$ into the four quadrants $B_0, B_1, B_2$ and $B_3$, where the boundaries of the $B_i$ are determined by $G^{-1}$. We assign random probability $\theta_i$ to quadrant $B_i$ with $(\theta_0, ..., \theta_3) \sim \textrm{Dir}(\alpha_0, ..., \alpha_3)$. Then we recursively partition $B_\kappa$ into quadrants $B_{\kappa 0}, ..., B_{\kappa 3}$, and split the mass assigned to $B_\kappa$ according to $(\theta_{\kappa 0}, ..., \theta_{\kappa 3})\sim\textrm{Dir}(\alpha_{\kappa 0}, ..., \alpha_{\kappa 3})$. This partitioning scheme is shown in Figure \ref{2d_partition_diagram_1}. We will denote this two-dimensional \textit{canonical} family of partitions with $\mathcal{T}_2$, the set of parameters $\alpha_\kappa$ with $\mathcal{A}_2$, and the set of splitting variables $\theta_\kappa$ with $\Theta_2$, where the subscript $_2$ emphasises the dimension of the space $\mathcal{X}\times \mathcal{Y}$. This leads to the following definition of the two-dimensional P\'olya tree:

\begin{definition}[\citet{hanson2006}]
    A random probability measure $\mathcal{P}$ on $(\mathcal{X}\times\mathcal{Y}, \mathcal{B}(\mathcal{X}\times\mathcal{Y}))$ is said to have a P\'olya tree distribution with parameter $(\mathcal{T}_2, \mathcal{A}_2)$, written $\mathcal{P} \sim {\normalfont\textrm{PT}}(\mathcal{T}_2, \mathcal{A}_2)$, if there exist nonnegative numbers $\mathcal{A}_2 = \{(\alpha_{\kappa 0}, \alpha_{\kappa 1}, \alpha_{\kappa 2}, \alpha_{\kappa 3}): \kappa \in K_2\}$ and random variables $\Theta_2 = \{(\theta_{\kappa 0}, \theta_{\kappa 1}, \theta_{\kappa 2}, \theta_{\kappa 3}) : \kappa \in K_2\}$ such that the following hold:
\begin{enumerate}[itemsep=1.8pt, parsep=1.8pt]
    \item all the random variables in $\Theta_2$ are independent;
    \item for every $\kappa\in K_2$ we have $(\theta_{\kappa 0}, \theta_{\kappa 1}, \theta_{\kappa 2}, \theta_{\kappa 3}) \sim {\normalfont\textrm{Dir}}(\alpha_{\kappa 0}, \alpha_{\kappa 1}, \alpha_{\kappa 2}, \alpha_{\kappa 3})$;
    \item for every $j \in \NN$ and every $\kappa \in \{0, 1, 2, 3\}^j$ we have $\mathcal{P}(B_{\kappa} | \Theta_2) = \prod_{i=1}^j \theta_{\kappa_1 ...\kappa_{i}}$.
\end{enumerate}
\end{definition}

\begin{figure}
    \centering
    \begin{subfigure}[t]{0.49\textwidth}
        \centering
        \fbox{\input{graphics/2d_partition_1}}
        \caption{Partitioning scheme for $X\dep Y$.}
        \label{2d_partition_diagram_1}
    \end{subfigure}
     \begin{subfigure}[t]{0.49\textwidth}
        \centering
        \fbox{\input{graphics/2d_partition_2}}
        \caption{Partitioning scheme for $X\indep Y$.}
        \label{2d_partition_diagram_2}
    \end{subfigure}
\end{figure}

Similarly to the one-dimensional case, samples from the P\'olya tree $\mathcal{P} \sim \textrm{PT}(\mathcal{T}_2, \mathcal{A}_2)$ are continuous with respect to the two-dimensional Lebesgue measure if we take $\alpha_{\kappa 0} = \alpha_{\kappa 1} = \alpha_{\kappa 2} = \alpha_{\kappa 3} = (j + 1)^2$, where $j$ denotes the length in the string $\kappa \in K_2$ \citep{walker1999}. Similar to the one-dimensional case, we only consider partitions up to a pre-specified depth $J(n)$.

When observing a sample $(X_1, Y_1), ..., (X_n, Y_n)$ from continuous random variables $X$ and $Y$ of which the joint distribution lies in the support of the two-dimensional P\'olya tree $\mathcal{P}$, we have that the marginal likelihood of that sample is 
\begin{equation}
    p((X, Y)_{1:n} | \Theta_2, \mathcal{T}_2, \mathcal{A}_2) = \prod_{\kappa\in K}\theta_{\kappa 0}^{n_{\kappa 0}}\theta_{\kappa 1}^{n_{\kappa 1}} \theta_{\kappa 2}^{n_{\kappa 2}} \theta_{\kappa 3}^{n_{\kappa 3}}.
\end{equation}
If we integrate over all possible values of all $\theta_\kappa$, we obtain the marginal likelihood

\begin{equation}
    p((X, Y)_{1:n} | \mathcal{T}_2, \mathcal{A}_2) = \prod_{\kappa \in K}\frac{\tilde{\mathrm{B}}(n_{\kappa 0} + \alpha_{\kappa 0}, n_{\kappa 1} + \alpha_{\kappa 1}, n_{\kappa 2} + \alpha_{\kappa 2}, n_{\kappa 3} + \alpha_{\kappa 3})}{\tilde{\mathrm{B}}(\alpha_{\kappa 0}, \alpha_{\kappa 1}, \alpha_{\kappa 2}, \alpha_{\kappa 3})},
\end{equation}

where $\tilde{\mathrm{B}}$ denotes the multivariate Beta function.\footnote{which is defined as $\tilde{\mathrm{B}}(\alpha_1, \alpha_2, \alpha_3, \alpha_4) := \left.\prod_{i=1}^4\Gamma(\alpha_i) \right/ \Gamma(\sum_{i=1}^4\alpha_i)$.}

Under the assumption $X\indep Y$, we construct a prior similar to the two-dimensional P\'olya tree. First we note that the two-dimensional family of partitions $\mathcal{T}_2$ can be regarded as the per-level Cartesian product of the partitions, i.e.\ 

\begin{equation}
    \mathcal{T}_2 = \big\{\{B_{\kappa} \times B_{\ell} : B_{\kappa} \in \mathcal{T}_X, B_{\ell} \in \mathcal{T}_Y, \kappa, \ell \in \{0,1\}^j\} : j\in\NN\big\}
\end{equation}

where $\mathcal{T}_X$ and $\mathcal{T}_Y$ are one-dimensional canonical partitions $\mathcal{X}$ and $\mathcal{Y}$ respectively. For every level $\kappa$, we first split the mass over the elements of $\mathcal{T}_X$ according to $(\theta_{\kappa 0}^X, \theta_{\kappa 1}^X)\sim \textrm{Dir}(\alpha_{\kappa 0}^X, \alpha_{\kappa 1}^X)$, and then independently split the mass over the elements of $\mathcal{T}_Y$ according to $(\theta_{\kappa 0}^Y, \theta_{\kappa 1}^Y)\sim \textrm{Dir}(\alpha_{\kappa 0}^Y, \alpha_{\kappa 1}^Y)$. We denote the set of parameters $\alpha_{\kappa}^X$ with $\mathcal{A}_X$, and the parameters $\alpha_{\kappa}^Y$ with $\mathcal{A}_Y$. This prior yields a marginal likelihood of
\begin{align}\label{2d_likelihood}
\begin{split}
    p((X, Y)_{1:n} | \mathcal{T}_2, \mathcal{A}_X, \mathcal{A}_Y) & = \prod_{\kappa \in K}\frac{\textrm{B}(n_{\kappa 0} + n_{\kappa 2} + \alpha_{\kappa 0}^X, n_{\kappa 1} + n_{\kappa 3} + \alpha_{\kappa 1}^X)}{\textrm{B}(\alpha_{\kappa 0}^X, \alpha_{\kappa 1}^X)} \\
    & \quad\quad\quad \times 
    \frac{\textrm{B}(n_{\kappa 0} + n_{\kappa 1} + \alpha_{\kappa 0}^Y, n_{\kappa 2} + n_{\kappa 3} + \alpha_{\kappa 1}^Y)}{\textrm{B}(\alpha_{\kappa 0}^Y, \alpha_{\kappa 1}^Y)},
\end{split}
\end{align}
as shown by \citet{filippi2017}. We notice that this equals the product of the marginal likelihoods of $X$ and $Y$ according to independent one-dimensional P\'olya tree priors $\mathcal{P}_X\sim \textrm{PT}(\mathcal{T}_X, \mathcal{A}_X)$ on $\mathcal{X}$ and $\mathcal{P}_Y\sim \textrm{PT}(\mathcal{T}_Y, \mathcal{A}_Y)$ on $\mathcal{Y}$, i.e.
\begin{equation}
\label{2d_likelihood_indep}
    p((X, Y)_{1:n} | \mathcal{T}_2, \mathcal{A}_X, \mathcal{A}_Y) = p(X_{1:n}|\mathcal{T}_X, \mathcal{A}_X)p(Y_{1:n}|\mathcal{T}_Y, \mathcal{A}_Y),
\end{equation}
where the univariate marginal likelihoods are computed according to equation (\ref{app:1d_likelihood}). To ensure that this prior is not biased when considered in conjunction with the two-dimensional P\'olya tree, we consider parameters $\alpha_{\kappa 0}^X = \alpha_{\kappa 0} + \alpha_{\kappa 2}$, $\alpha_{\kappa 1}^X = \alpha_{\kappa 1} + \alpha_{\kappa 3}$, $\alpha_{\kappa 0}^Y = \alpha_{\kappa 0} + \alpha_{\kappa 1}$ and $\alpha_{\kappa 1}^Y = \alpha_{\kappa 2} + \alpha_{\kappa 3}$ \citep{filippi2017}. When using the set of standard parameters $\mathcal{A}_2$ for the two-dimensional P\'olya tree, we have $\mathcal{A}' := \mathcal{A}_X = \mathcal{A}_Y = \{2j^2 : j\in \NN\}$.

\subsection{A nonparametric independence test}
\label{independence_test}
A Bayesian independence test that utilises two-dimensional P\'olya trees is proposed by \citet{filippi2017}. Considering one-dimensional continuous random variables $X$ and $Y$, we test the hypotheses
\begin{equation}
    H_0: X \indep Y, \quad H_1: X \dep Y
\end{equation}
using the Bayes factor
\begin{equation}
    \mathrm{BF}_{01} = \frac{p(X_{1:n} | \mathcal{T}, \mathcal{A}')p(Y_{1:n} | \mathcal{T}, \mathcal{A}')}{p((X, Y)_{1:n} | \mathcal{T}_2, \mathcal{A}_2)},
\end{equation}
where the marginal likelihoods are computed according to equations (\ref{app:1d_likelihood}) and (\ref{2d_likelihood}).

Using similar arguments as for the two-sample test, the Bayes factor can be denoted as an infinite product, of which the terms are equal to one when $n_{XY|\kappa} \leq 1$. Therefore we compute the marginal likelihoods up to level $J(n)$, or until all elements of the partition contain at most one observation.

%% file: graphics/2d_partition_1.tex
\begin{tikzpicture}[xscale=0.8, yscale=0.8]

    \node () at (0, 5.2) {};
  
    \draw (-2, 5) -- (2, 5) ;
    \draw (-2, 3) -- (2, 3) ;
    \draw (-2, 1) -- (2, 1) ;
    \draw (-2, 1) -- (-2, 5) ;
    \draw (0, 1) -- (0, 5) ;
    \draw (2, 1) -- (2, 5) ;
  
    \node () at (-2.5, 3) {$\mathcal{Y}$};
    \node () at (0, 0.65) {$\mathcal{X}$};
  
    \node () at (-1.5, 1.4) {\footnotesize $B_0$};
    \node () at (1.5, 1.4) {\footnotesize $B_1$};
    \node () at (-1.5, 4.6) {\footnotesize $B_2$};
    \node () at (1.5, 4.6) {\footnotesize $B_3$};
  
    \draw[->, line width=0.8pt] (0, 3) -- (-1, 2) node [pos=0.5, left=0.5em] {\footnotesize $\theta_{0}$};
    \draw[->, line width=0.8pt] (0, 3) -- (1, 2) node [pos=0.5, right=0.5em] {\footnotesize $\theta_{1}$};
    \draw[->, line width=0.8pt] (0, 3) -- (-1, 4) node [pos=0.5, left=0.5em] {\footnotesize $\theta_{2}$};
    \draw[->, line width=0.8pt] (0, 3) -- (1, 4) node [pos=0.5, right=0.5em] {\footnotesize $\theta_{3}$};
  
    \node () at (3.5, 3) [text width=1.4cm, align=left] {\footnotesize Level 1 \\ for $X \dep Y$};

  \end{tikzpicture}

%% file: graphics/2d_partition_2.tex
\begin{tikzpicture}[xscale=0.8, yscale=0.8]
    \node () at (0, 5.2) {};
  
    \draw (-2, 5) -- (2, 5) ;
    \draw (-2, 3) -- (2, 3) ;
    \draw (-2, 1) -- (2, 1) ;
    \draw (-2, 1) -- (-2, 5) ;
    \draw (0, 1) -- (0, 5) ;
    \draw (2, 1) -- (2, 5) ;
  
    \node () at (-2.5, 3) {$\mathcal{Y}$};
    \node () at (0, 0.65) {$\mathcal{X}$};
  
    \node () at (-1.5, 1.4) {\footnotesize $B_0$};
    \node () at (1.5, 1.4) {\footnotesize $B_1$};
    \node () at (-1.5, 4.6) {\footnotesize $B_2$};
    \node () at (1.5, 4.6) {\footnotesize $B_3$};
  
    \draw[->, line width=0.8pt] (0, 3) -- (-1, 2) node [pos=0.5, left=0.5em] {\footnotesize$\theta^{X}_{0}\theta^{Y}_{0}$};
    \draw[->, line width=0.8pt] (0, 3) -- (1, 2) node [pos=0.5, right=0.5em] {\footnotesize$\theta^{X}_{1}\theta^{Y}_{0}$};
    \draw[->, line width=0.8pt] (0, 3) -- (-1, 4) node [pos=0.5, left=0.5em] {\footnotesize$\theta^{X}_{0}\theta^{Y}_{1}$};
    \draw[->, line width=0.8pt] (0, 3) -- (1, 4) node [pos=0.5, right=0.5em] {\footnotesize$\theta^{X}_{1}\theta^{Y}_{1}$};
  
    \node () at (3.5, 3) [text width=1.4cm, align=left] {\footnotesize Level 1 \\ for $X \indep Y$};
  
  \end{tikzpicture}

%% file: input/B_sachs.tex
\begin{figure}[ht]
    \centering
    \includegraphics[width=0.7\textwidth]{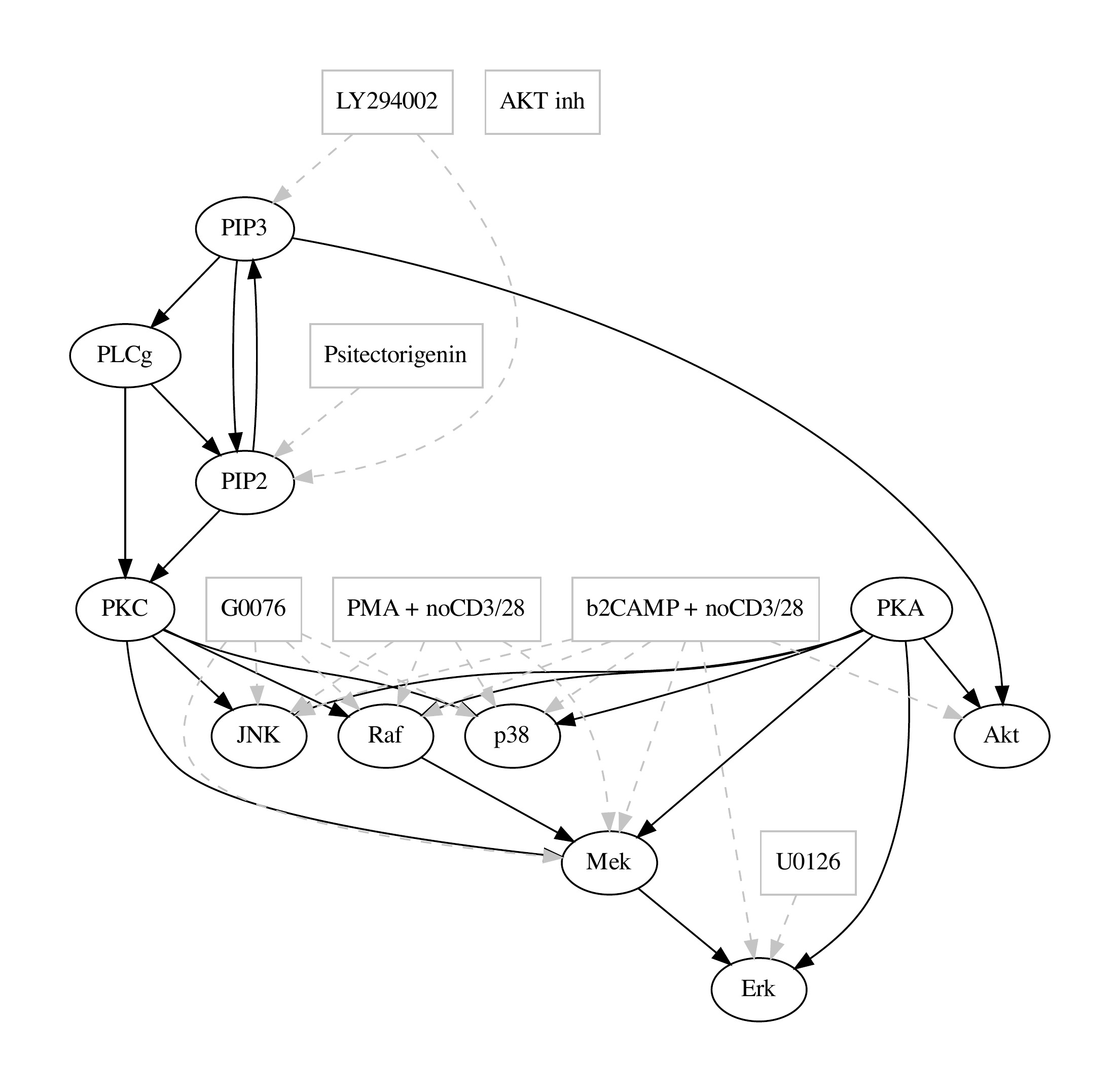}
    \caption{The ‘expert network’ as provided by \citet{sachs2005}. Edges indicate direct causal effects between the nodes. Interventions and their direct causal effects are indicated with light-coloured and dashed nodes and edges.}
    \label{expert_network}
\end{figure}

In the main paper we apply the LCD algorithm, implemented with the Bayesian ensemble of independence tests, to protein expression data \citep{sachs2005}. The data set consists of measurements of 11 phosphorylated proteins and phospholipids (Raf, Erk, p38, JNK, Akt, Mek, PKA, PLCg, PKC, PIP2 and PIP3) and 8 indicators of different interventions, performed by adding reagents to the cellular system, which are depicted in Table \ref{interventions}.\footnote{Similarly to most analyses of this data, we restrict our attention to 8 out of 14 experimental conditions, namely those in which no ICAM was added.} The biological details of these proteins, phospholipids, and reagents are described in \citet{sachs2005}. Using flow cytometry, the activity of the 11 proteins and phospholipids are measured from a single human immune system cell. Flow cytometry allows for simultaneous, independent observation of hundreds of cells, producing a statistically large sample, and thus allowing for the application of causal inference algorithms \citep{sachs2005}. The ‘expert network’ from \citet{sachs2005} is depicted in Figure \ref{expert_network}. We note that, as argued in the main paper, we do not accept this network as the true causal graph, but merely display it suggestively.

\begin{table}[h]
    \caption{Interventions from the data set of \citet{sachs2005}.}
    \label{interventions}
    \centering
    \begin{tabular}{llr}
        \toprule
         & Description & Nr. of observations \\
        \midrule
        1 & CD3, CD28 & 853 \\
        2 & CD3, CD28, Akt-inhibitor & 911 \\
        3 & CD3, CD28, G0076 & 723 \\
        4 & CD3, CD28, Psitectorigenin & 810 \\
        5 & CD3, CD28, U0126 & 799 \\
        6 & CD3, CD28, LY294002 & 848 \\
        7 & PMA & 913 \\
        8 & $\beta$2CAMP & 707 \\
        \bottomrule
    \end{tabular}
\end{table}

We assume that adding the reagents is not caused by the activity of the proteins and phospholipids, which justifies the application of the LCD algorithm to this dataset, as per Proposition 3.1 of the main paper. When performing a statistical test we always use the entire set of observations. As is common when analysing flow cytometry data, we preprocessed the data by taking the log of the raw measurement values.
\clearpage

